\documentclass[twocolumn]{autart}
\usepackage{amssymb,amsmath,dsfont}
\usepackage{graphicx}
\usepackage{subfigure,xcolor}
\usepackage{float}

\usepackage[english]{babel}
\usepackage{caption}
\usepackage{multirow}

\newcommand {\etavec}{\boldsymbol{\eta}}

\newcommand {\phivec}{{\boldsymbol{\phi}}}
\newcommand {\Phivec}{{\boldsymbol{\Phi}}}
\newcommand {\varphivec}{{\boldsymbol{\varphi}}}
\newcommand {\betavec}{{\boldsymbol{\beta}}}

\newcommand {\col}{{\mathrm{col}}}

\newcommand {\alphavec}{\boldsymbol{\alpha}}

\newcommand {\xivec}{{\boldsymbol{\xi}}}
\newcommand {\lambdavec}{{\boldsymbol{\lambda}}}

\newcommand {\zetavec}{{\boldsymbol{\zeta}}}

\newcommand {\Psivec}{\mbox{\boldmath $\Psi$}}

\newcommand {\thetavec}{{\boldsymbol{\theta}}}

\newfont{\pseudocode}{cmtt10}

\newcommand{\bfx}{\mathbf{x}}
\newcommand{\bfl}{\boldsymbol{\ell}\,}

\newcommand{\bff}{\mathbf{f}}
\newcommand{\bfv}{\mathbf{v}}
\newcommand{\bfb}{\mathbf{b}}
\newcommand{\bfB}{\mathbf{B}}
\newcommand{\bfJ}{\mathbf{J}}
\newcommand{\bfe}{\mathbf{e}}
\newcommand{\bfd}{\mathbf{d}}
\newcommand{\bfA}{\mathbf{A}}
\newcommand{\bfI}{\mathbf{I}}
\newcommand{\bfR}{\mathbf{R}}

\newcommand{\bfM}{\mathbf{M}}

\newcommand{\bfC}{\mathbf{C}}
\newcommand{\bfz}{\mathbf{z}}

\newcommand{\bfh}{\mathbf{h}}

\newcommand{\bfs}{\mathbf{s}}
\newcommand{\bfu}{\mathbf{u}}
\newcommand{\bfa}{\mathbf{a}}
\newcommand{\bfg}{\mathbf{g}}

\newcommand{\bfp}{\mathbf{p}}
\newcommand{\bfq}{\mathbf{q}}
\newcommand{\bfP}{\mathbf{P}}
\newcommand{\bfQ}{\mathbf{Q}}
\newcommand{\bfG}{\mathbf{G}}

\newcommand{\dist}{\mathrm{dist}}

\newcommand{\pd}{{\partial}}

\newcommand{\Real}{\mathbb{R}}
\newcommand{\Natural}{\mathbb{N}}
\newcommand{\Numbers}{\mathbb{Z}}

\newcommand{\norm}[1]{\left\Vert#1\right\Vert}

\newcommand {\ess}{{\mathrm{ess}}}

\renewcommand{\Real}{\mathbb{R}}

\renewcommand{\Real}{\mathds{R}}
\renewcommand{\Natural}{\mathds{N}}

\renewcommand{\Numbers}{\mathds{Z}}

\newcommand{\ra}{\rightarrow}

\newcommand{\Ecal}{\mathcal{E}}

\renewcommand{\col}[1]{\mathrm{col}\left(#1\right)}

\newcommand{\beq}{\begin{equation}}
\newcommand{\eeq}{\end{equation}}

\newcommand{\splt}[1]{\begin{split}#1\end{split}}

\newcommand{\bbm}{\begin{bmatrix}}
\newcommand{\ebm}{\end{bmatrix}}

\newcommand{\bpm}{\begin{pmatrix}}
\newcommand{\epm}{\end{pmatrix}}

\newcommand{\bit}{\begin{itemize}}
\newcommand{\eit}{\end{itemize}}

\newcommand{\ben}{\begin{enumerate}}
\newcommand{\een}{\end{enumerate}}

\newcommand{\barr}{\begin{array}}
\newcommand{\earr}{\end{array}}

\renewcommand{\norm}[1]{\Vert#1\Vert}

\newtheorem{assume}{Assumption}[section]

\newenvironment{proof}[1][Proof]{\begin{trivlist}
\item[\hskip \labelsep {\bfseries #1}]}{\end{trivlist}}

\begin{document}

\begin{frontmatter}

\title{Adaptive Observers and Parameter Estimation for a Class\\  of Systems Nonlinear in the Parameters}

\thanks[footnoteinfo]{This paper was not presented at any IFAC
meeting.  Cees van Leeuwen is supported by an Odysseus grant from the Flemish Science Organization FWO. Corresponding author I.~Yu.~Tyukin. Tel.
+44-116-2525106.}

\author[First]{Ivan Y. Tyukin}
\author[Second]{Erik Steur}
\author[Third]{Henk Nijmeijer}
\author[Second]{Cees van Leeuwen}

\address[First]{Dept. of Mathematics, University of Leicester, Leicester, LE1 7RH, UK (Tel: +44-116-252-5106; I.Tyukin@le.ac.uk)}
\address[Second]{Laboratory for Perceptual Dynamics, KU Leuven, Tiensestraat 102, 3000 Leuven, Belgium (erik.steur@ppw.kuleuven.be,  cees.vanleeuwen@ppw.kuleuven.be)}
\address[Third]{Dept. of Mechanical Engineering, Eindhoven University of Technology, P.O. Box 513 5600 MB,  Eindhoven, The Netherlands, (h.nijmeijer@tue.nl)}



\begin{abstract}
    We consider the problem of asymptotic reconstruction of the
    state and parameter values in systems of ordinary differential equations. A
    solution to this problem is proposed for a class of systems
    of which the unknowns are allowed to be nonlinearly
    parameterized functions of state and time.  Reconstruction of state and
    parameter values is  based on the concepts
    of weakly attracting sets and non-uniform convergence and is subjected to persistency of excitation conditions.
    In the absence of nonlinear parametrization the resulting
    observers reduce to standard estimation schemes. In this respect, the proposed method constitutes a generalization
    of the conventional canonical adaptive observer design.
\end{abstract}

\end{frontmatter}


\section{Introduction}

We consider observer-based methods for state and parameter
estimation in nonlinear dynamical systems. These methods are
effective as long as the original system has, or can be
transformed into, one of the {\it canonical adaptive observer
forms} \cite{Bastin88}, \cite{Marino90}, \cite{Besancon:2000}.
Their common characteristic is linearity in the unknown parameters.  For this class of
systems, subject to persistency of excitation conditions, reconstruction of state and parameter vectors can be achieved exponentially fast.
%
%

There are systems, however, in which the unknown parameters
enter the model nonlinearly. These systems constitute a
remarkably wide class including models in
chemical kinetics \cite{Gorban:2005}, \cite{Bastin:1990},
biology and neuroscience \cite{Izhikevich:2007}. Whereas the
problem of state estimation can be solved for a large class of
nonlinearly parameterized systems \cite{MarinoTomei93_2},
observer-based parameter reconstruction
 is often confined to systems with  monotone
\cite{IEEE_TAC:2007:tyukin},  \cite{Ortega2011} or one-to-one
parameterizations \cite{IEEE:2009:Grip},
\cite{Automatica:Grip:2010}, \cite{SysConLett:Grip:2011},
\cite{Automatica:Farza:2009}.

Several authors have recently advanced strategies for
overcoming these limitations. In
\cite{Automatica:2008:Johnson},  combining interval
analysis  with multiple shooting methods is proposed
to tackle the state and  parameter reconstruction problem.
Another interesting approach is presented in
\cite{Abarbanel:2009}: the original continuous-time model is
replaced with a discrete-time approximation. Measured variables
are then considered as known functions of unknown parameters
and initial conditions, of which the estimates can be found by off-line nonlinear optimization routines (see
also \cite{Rao:2003}, \cite{Alessandri:2008} where optimization
techniques with moving horizon are discussed). These
approaches offer obvious advantages, e.g.
the availability of a vast library of numerical methods for
solving general nonlinear  optimization problems. Nevertheless,
these methods run into restrictions too.  Exhaustive search for
a global minimum may become intractable for dimensions higher
than $1$ or $2$. On the other hand, if conventional
polynomial-complexity algorithms are used then the possibility arises that the algorithm will converge to a local minimum.

In this paper we explore further  possibilities of developing
adaptive observers for systems which are both linearly and
nonlinearly parameterized.
The parametrization is not required to be invertible or
monotone. Our approach combines the advantages of the existing
schemes, in being capable of ensuring exponentially fast
convergence with the flexibility of explorative behavior, a
behavior inherent to algorithms for solving genuine nonlinear
optimization problems. Inference of
the values of state and a part of the parameter vector of these
systems is achieved by employing exponentially fast
converging estimators. Estimation of the values of the
remaining parameters is based on an explorative search procedure.
%
Since exploration is restricted to a subset of
the unknown parameters, the proposed strategy reduces the
overall computational costs, as compared to when full-scale search-based optimization had been invoked.

The resulting observer can be imagined as a system comprising
of an exponentially stable part coupled with an explorative
one. Systems of this type  have previously been used
in adaptive control
\cite{CDC:Pomet:1992}, \cite{Dyn_Con:Ilchman:97},
\cite{SysConLett:Mertensson:85},
\cite{SysConLett:Mertensson:93}.
Here we
demonstrate that
these classical ideas can be
applied
to the problem of adaptive
observer design for systems which are nonlinearly dependent on
parameters. We show that, subject to a condition of persistent
excitation, it is possible to reconstruct state and parameters
of a reasonably broad subclass of these systems.

The paper is organized as follows. Notational agreements are
introduced in Section \ref{sec:notation}. Section
\ref{sec:problem} provides the formal statement of the problem,
Sections \ref{sec:results}, \ref{sec:proof} contain main
results of the article, Section \ref{sec:discussion} discusses
possible generalizations, Section \ref{sec:example} contains illustrative
examples, and Section \ref{sec:conclusion} concludes the paper.
Proofs of auxiliary results are presented in the Appendix.


\section{Notation}\label{sec:notation}

The following notational conventions are used throughout the
paper: \bit
    \item $\Real$ denotes the set of real numbers, $\Real_{>a} = \{x\in\Real\ |\ x>a\}$, and  $\Real_{\geq a} = \{x\in\Real\ |\ x\geq a\}$;
    \item $\Numbers$ denotes the set of integers, and $\Natural$ stands for the set of positive integers;
    \item the Euclidean norm of $\bfx \in \Real^n$ is denoted by $\norm{\bfx}$, $\norm{\bfx}^2 = \bfx^T \bfx$, where ${}^T$ stands for transposition;
      \item the space of $n\times n$ matrices with real entries is denoted by $\Real^{n\times n}$; let  $\bfP\in\Real^{n\times n}$, then $\bfP>0$ ($\bfP\geq 0$) indicates that $\bfP$ is symmetric and positive (semi-)definite; $\bfI_{n}$ denotes  the $n\times n$ identity matrix.
    \item  by ${L}^n_\infty[t_0,T]$, $t_0\in\Real$, $T\in\Real, \ T\geq t_0$ we denote the space of all functions $\bff:[t_0,T]\rightarrow\Real^n$ such that $\|\bff\|_{\infty,[t_0,T]}=\ess \sup\{\|\bff(t)\|,t \in [t_0,T]\}<\infty$; $\|\bff\|_{\infty,[t_0,T]}$ stands for the ${L}^n_\infty[t_0,T]$ norm of $\bff(t)$; if the function $\bff$ is defined  on a set larger than $[t_0,T]$ then notation $\|\bff\|_{\infty,[t_0,T]}$ applies to the restriction of $\bff$ on $[t_0,T]$;
    \item $\mathcal{C}^{r}$ denotes the space of continuous functions that are at least $r$ times differentiable;
    \item Let $\mathcal{A}$ be a subset of $\Real^n$, then for
        all $\bfx\in\Real^n$, we define
        $\mathrm{dist}(\mathcal{A},\bfx)=\inf_{\bfq\in\mathcal{A}}\|\bfx-\bfq\|$;
        { \item A solution of
        $\dot{\bfx}=\bff(t,\bfx,\thetavec,u(t))$,
        $\bff:\Real\times\Real^n\times\Real^m\times\Real\rightarrow\Real^n$, $\thetavec\in\Real^m$,
        $u:\Real\rightarrow\Real$ passing through
        $\bfx_0\in\Real^n$ at $t=t_0$ is denoted by
        $\bfx(t,t_0,\bfx_0,\thetavec,[u])$. In cases when
        $u$ and/or $\thetavec$, $\bfx_0$, $t_0$ are clearly determined by the context a more compact
        notation, $\bfx(t,t_0,\bfx_0,\thetavec)$ (or $\bfx(t,t_0,\bfx_0)$, $\bfx(t)$, $\bfx$ respectively), is
                used.}
    \item The symbol $\mathcal{K}$ denotes the class of all strictly increasing continuous functions $\kappa:\Real_{\geq 0}\rightarrow\Real_{\geq 0}$ such that $\kappa(0)=0$; the symbol $\mathcal{K}_{\infty}$ denotes the class of all functions  $\kappa\in\mathcal{K}$ such that $\lim_{s\rightarrow\infty}\kappa(s)=\infty$.
    \item Let $\epsilon \in \Real_{\geq 0}$, then $\norm{\bfx}_\epsilon$ stands for: $\norm{\bfx}-\epsilon$ if $\norm{\bfx} > \epsilon$, and $0$ otherwise.
    \item Finally, for $\lambdavec\in\Real^p$ and $\thetavec\in\Real^m$, the notation $(\lambdavec,\thetavec)$ stands for  $\mathrm{col}(\lambda_1,\ldots,\lambda_p,\theta_1,\ldots,\theta_m)$.
\eit


\section{Preliminaries and Problem formulation}\label{sec:problem}
\subsection{Adaptive observer canonical form}

Throughout the paper we will focus exclusively on the class of
systems
that are forward-complete: 
\begin{defn}
Let $\mathcal{L}_u[t_0,\infty]$ be a subspace of
$L_{\infty}[t_0,\infty]$. A single-input single-output system
described by $\dot{\bfx}=\bff(t,\bfx,u(t)),\
\bff:\Real\times\Real^n\times\Real\rightarrow\Real^n$,
$y=\bfh(t,\bfx),  \ \bfh:\Real\times\Real^n\rightarrow\Real, \
u\in \mathcal{L}_u[t_0,\infty]$, where $u$ is the input, and
$y$ is the output, is called forward-complete (with respect to
$\mathcal{L}_u$) iff for  any $t_0\in\Real$, $\bfx_0\in\Real^n$, and
$u\in\mathcal{L}_u[t_0,\infty]$ the solution
$\bfx(t,t_0,\bfx_0,[u])$ exists and is defined for all $t\geq
t_0$.
\end{defn}
 Let
$\mathcal{L}_{u}[t_0,\infty]=L_{\infty}[t_0,\infty]\cap
\mathcal{C}^{0}[t_0,\infty]$, and consider a forward-complete
single-input single-output system;  let $\bfx\in\Real^n$ be its state,
$y:\Real\rightarrow\Real$ be the measured output, and
$u:\Real\rightarrow\Real$, $u\in\mathcal{L}_u$ be the input. We
recall that a system is in the {\it adaptive observer canonical
form} if it is governed by the
following set of equations{
\begin{equation}\label{eq:observer_canonical_form}
    \splt{
        \dot{\bfx} &= \bfA \bfx + \bfB\phivec^T(t,y) \thetavec + \bfg(t,y,u(t)) \\
        y &= \bfC^T \bfx, \ \bfx(t_0)=\bfx_0, \ \bfx_0\in\Real^n,
    }
\end{equation}
where
$\bfA=\left(\begin{array}{c|c}
\multirow{2}{*}{\ensuremath{\bfa}} &   I_{n-1}\\ & 0 \end{array}\right)$, $\bfa\in\Real^n$, $\bfB\in\Real^n$,  $\bfC\in\Real^n$,
$\bfB = \mathrm{col}(1, b_1,\dots,b_{n-1})$,  $\bfC =
\mathrm{col}(1, 0, \dots,0)$, } the functions
$\phivec:\Real\times\Real\rightarrow\Real^m$,
$\bfg:\Real\times\Real\times\Real\rightarrow \Real^n$,
$\phivec,\bfg\in\mathcal{C}^0$ are known, and $\thetavec\in\Real^m$ is the vector of
the unknown parameters.  The triplet $\bfA,\bfB,\bfC$ is
supposed to satisfy
\begin{equation}\label{eq:condition_MKY}
         \left\{\begin{split}
         &\bfP(\bfA+\bfl \bfC^T)+(\bfA+\bfl \bfC^T)^T \bfP\leq - \bfQ \\
         &\bfP\bfB = \bfC.
         \end{split}\right.
\end{equation}
for some  $\bfl\in\Real^n$ and  $\bfP>0$, $\bfQ>0$. Although
condition (\ref{eq:condition_MKY}) may appear restrictive,
it has been shown in \cite{Marino92} that subject
to the very natural constraint that the pair $\bfA,\bfC$ is
observable, there is
a time-varying
parameter-dependent coordinate transformation such that in new
coordinates the system is still of the form
(\ref{eq:observer_canonical_form}) and satisfies condition
(\ref{eq:condition_MKY}). If requirement
(\ref{eq:condition_MKY}) holds then the system
\begin{equation}\label{eq:canonic_observer}
\begin{split}
\dot{\hat{\bfx}}=&\bfA\hat{\bfx}+\bfl(\bfC^T\hat{\bfx}-y(t))+\bfB\phivec^T(t,y(t)) \hat{\thetavec}\\
& + \bfg(t,y(t),u(t))\\
\dot{\hat{\thetavec}}=&-\gamma(\bfC^T\hat{\bfx}-y(t)) \phivec(t,y(t)), \ \gamma\in\Real_{>0},
\end{split}
\end{equation}
where $\hat \bfx(t)\in\Real^n$, $\hat \thetavec(t)\in\Real^m$,
is an adaptive observer for (\ref{eq:observer_canonical_form})
(cf. \cite{Marino90}, \cite{Bastin88}) provided that the restriction of
$\phivec(\cdot,y(\cdot))$ on $\Real_{\geq t_0}$ is {\it
persistently exciting}:
\begin{defn}\label{defn:pe} A function $\betavec:\ {\Real_{\geq t_0}}  \ra \Real^{m}$ is said to be persistently exciting if there exist
$L,\mu\in\Real_{>0}$:
    \beq\label{eq:PE_linear_uniform:1}
       \begin{array}{l} \int_{t}^{t+L} \betavec(\tau)\betavec^T(\tau){\rm d}\tau \geq \mu
        \bfI_m, \ \forall \ t\geq t_0.\end{array}
    \eeq
\end{defn}
The fact that (\ref{eq:canonic_observer}) is an adaptive
observer for (\ref{eq:observer_canonical_form}) is based on a
well-known result on the exponential stability of the following
class of linear time-varying systems\footnote{In the context of
adaptive observer design for (\ref{eq:observer_canonical_form})
the function $\betavec$ in (\ref{eq:error_dynamics}) is
defined as $\betavec(t)=\phivec(t,y(t))$, $t\geq t_0$.}
\begin{equation}\label{eq:error_dynamics}
\dot{\bfe}=\bfA(t)\bfe, \ \bfA(t)=\left(\begin{array}{cc} \bfA+\bfl\bfC^T & \bfB\betavec^T(t)\\ -\betavec(t)\bfC^T & 0 \end{array}\right).
\end{equation}
The result is provided in Theorem
\ref{thm:preliminary:exponential_skew_symmetric} below (see
e.g. \cite{Automatica:Loria:2004} for a proof).

\begin{thm}\label{thm:preliminary:exponential_skew_symmetric}
Consider system (\ref{eq:error_dynamics}). Suppose that
condition (\ref{eq:condition_MKY}) holds for certain
$\bfl\in\Real^n$, $\bfP>0$, $\bfQ>0$, the function
$\betavec(t)$ is persistently exciting, and
\begin{equation}\label{eq:derivative_bound}
\exists \ M\in\Real_{>0}: \ \ \max\{ \|\betavec(t)\|,\|\dot{\betavec}(t)\|\}\leq M \ \forall \ t\geq t_0.
\end{equation}
Let $\Phi(t,t_0)$, $\Phi(t_0,t_0)=I_{n+m}$, be the fundamental
solution matrix of (\ref{eq:error_dynamics}).
Then there exist $\rho,D\in\Real_{>0}$ such that
$\|\Phi(t_2,t_1)\bfp\|\leq  D e^{-\rho (t_2-t_1)}\|\bfp\|$ for all
$t_2\geq t_1\geq t_0$ and $\bfp\in\Real^{n+m}$.
\end{thm}
 The parameters $\rho$ and $D$
can be expressed explicitly as functions
of $M$, $\mu$, $L$, and  $\bfA$, $\bfB$, $\bfC$, $\bfl$ \cite{Automatica:Loria:2004}. By letting
$\bfe=\mathrm{col}(\hat{\bfx}-\bfx,\hat{\thetavec}-\thetavec)$
and taking (\ref{eq:observer_canonical_form}),
(\ref{eq:canonic_observer}) into account one can confirm that
the system-observer equations are of form
(\ref{eq:error_dynamics}). Thus, subject to persistency of
excitation of the restriction of $\phivec^{T}(\cdot,y(\cdot))$ on $\Real_{\geq t_0}$,
$\lim_{t\rightarrow\infty}\hat{\bfx}(t,t_0,\hat{\bfx}_0,\hat{\thetavec}_0)-{\bfx}(t,t_0,{\bfx}_0,\thetavec)=0$,
$\lim_{t\rightarrow\infty}\hat{\thetavec}(t,t_0,\hat{\bfx}_0,\hat{\thetavec}_0)=\thetavec$
along the solutions of (\ref{eq:observer_canonical_form}),
(\ref{eq:canonic_observer}), and the convergence is
exponential. The problem, however, is that { if some
parameters enter the equations nonlinearly then this}
%
creates an obstacle for the
explicit use of Theorem
\ref{thm:preliminary:exponential_skew_symmetric} and,
consequently, observer (\ref{eq:canonic_observer}). In the next
sections we present and analyze a class of systems nonlinear in the parameters which can be thought of as an immediate
generalization of (\ref{eq:observer_canonical_form}).

\subsection{Systems considered in this article}

We begin with the following class of forward-complete
single-input-single-output nonlinear systems:
\begin{equation}\label{eq:neural_model}
    \splt{
        \dot{\bfx} &= \bfA \bfx + \bfB\varphivec^T(t,\lambdavec,y) \thetavec + \bfg(t,\lambdavec,y,u(t)) + \xivec(t), \\
        y &= \bfC^T \bfx, \ \bfx(t_0)=\bfx_0, \ \bfx_0\in\Real^n,
        }
\end{equation}
where $\bfA\in\Real^{n\times n}$, and $\bfB, \bfC\in\Real^n$
are defined as in (\ref{eq:observer_canonical_form});
$\varphivec:\Real\times\Real^p\times\Real\rightarrow\Real^{m}$,
$\bfg:\Real\times\Real^p\times\Real\times\Real\rightarrow\Real^{n}$,
{are known continuous functions},
$\lambdavec=\mathrm{col}(\lambda_1,\dots,\lambda_p)\in\Real^p$,
$\thetavec=\mathrm{col}(\theta_1,\dots,\theta_m)\in\Real^{m}$
are {\it unknown parameters}, and
$u\in\mathcal{L}_u\cap \mathcal{C}^1$, is the input. We assume
that the values of $\lambdavec$, $\thetavec$ belong to the
hypercubes $\Omega_\lambda\subset\Real^p$,
$\Omega_\theta\subset\Real^{m}$ with known bounds:
$\theta_{i}\in[\theta_{i,\min},\theta_{i,\max}]$,
$\lambda_j\in[\lambda_{j,\min},\lambda_{j,\max}]$,
and that $y(t)\in\mathcal{D}_y$,
$u(t)\in\mathcal{D}_u$,
$\mathcal{D}_y,\mathcal{D}_u\subset\Real$ for $t\geq t_0$.

In \eqref{eq:neural_model},
$\bfx=\col{x_1,x_2,\dots,x_n}\in\Real^n$ is the state vector,
$y$ is the measured output, the input
$u$ is a known function, and
$\xivec\in\mathcal{C}^0:\Real\rightarrow\Real^n$ is an unknown
yet bounded continuous function:
\beq\label{eq:xi_bound}
  \exists \ \Delta_\xi\in\Real_{\geq 0}: \ \   \|\xivec(t)\|\leq \Delta_\xi \ \forall \ t,
\eeq
representing some {\it unmodeled dynamics} (e.g. noise).
The system's state $\bfx$ is not measured;  only the values of
the input $u(t)$ and the output $y(t)=x_1(t)$, $t\geq t_0$ in
(\ref{eq:neural_model}) are accessible over any time interval
$[t_0,t]$ that belongs to the history of the system.
\begin{table*}[!th]
\caption{Examples of physical systems in the form of (\ref{eq:neural_model}), (\ref{eq:neural_model:ext})}\label{tab:examples}
\centering
\begin{tabular}{|c|c|c|}
    \hline
       \begin{minipage}[h]{0.14\linewidth}
       \begin{center}
     Physical system
     \end{center}
     \end{minipage}
      &  \begin{minipage}[h]{0.32\linewidth}      \begin{center} Model \end{center}  \end{minipage} &
      \begin{minipage}[h]{0.46\linewidth}
      \begin{center}
      Possible parametrization: $\thetavec$, $\lambdavec$, $\bfA$, $\bfB$, $\varphivec$, $\Psivec$, $\bfg$
      \end{center}
      \end{minipage}\\
     \hline
   \begin{minipage}[h]{0.14\linewidth}
    1. Reactant dynamics in a reactor \cite{Comp:Poyton:2006} \end{minipage} & \begin{minipage}[h]{0.32\linewidth} \vskip -2mm
     \[\begin{array}{l}\dot{x}=\frac{F}{V}(x_0-x) -k_{\mathrm{ref}} e^{-\frac{E}{R}\left(\frac{1}{T(t)}-\frac{1}{T_{\mathrm{ref}}}\right)}x,\\ y= x \end{array}\]
     \end{minipage} &
     \begin{minipage}[h]{0.46\linewidth}      \vskip -2mm
     \vskip 2mm
      System (\ref{eq:neural_model}): $\thetavec=\col{\frac{F}{V}x_0, -\frac{F}{V},-k_{\mathrm{ref}}e^{\frac{E}{R T_{\mathrm{ref}}}}}, \ \lambda=-\frac{E}{R}$, $\bfA=0, \ \bfB=1, \ \varphivec(t,\lambda,y)=\col{1,y,e^{\lambda \frac{1}{T(t)}}y}, \ g=0$
     \end{minipage}   \\
    \hline
     \begin{minipage}[h]{0.14\linewidth} 2. Magnetic bearings \cite{ACC:2000:Lin} \end{minipage}&
     \begin{minipage}[h]{0.32\linewidth} \vskip -2mm
     \[\begin{array}{l} \ddot{x}+J^{-1} d(t)= b\left(\frac{q_2^2(I_2,a,x)}{(x+a)^2}-\frac{q_1^2(I_1,a,x)}{(x-a)^2}\right)\\ y=x,\ |x|<a-\varepsilon, \ a,\varepsilon\in\Real_{>0}, \end{array} \]
     \vskip -3mm $d,\dot{d}$ are bounded
     \end{minipage} &   \begin{minipage}[h]{0.46\linewidth} \vskip -2mm
        \vskip 2mm
         System (\ref{eq:neural_model}):
        $\thetavec=\col{-J^{-1}, b}, \ \bfA=\left(\begin{array}{cc} 0 & 1\\ 0 & 0\end{array}\right)$,  $\bfB=\left(\begin{array}{c}0\\1\end{array}\right)$, $\lambda=a$,
          $\varphivec(t,\lambda,y)=\col{d(t),\frac{q_2^2(I_2,\lambda,y)}{(y+\lambda)^2}-\frac{q_1^2(I_1,\lambda,y)}{(y-\lambda)^2}}, \ g=0$
     \end{minipage}  \\
    \hline
 \begin{minipage}[h]{0.14\linewidth} 3. Action potentials in a cell \cite{Rowat} \end{minipage}&
     \begin{minipage}[h]{0.32\linewidth} \vskip -5mm
     \[
     \begin{array}{l}
    \dot{x}_1=-\frac{x_1}{\tau_m} + \frac{A_f}{\tau_m}\tanh\left(\frac{\sigma_f}{A_f}x_1\right)-\frac{x_2}{\tau_m}\\
     \dot{x}_2=-\frac{x_2}{\tau_s}+\frac{\sigma_s}{\tau_s} x_1\\
     y=x_1
     \end{array}
     \]
     \end{minipage} &
  \begin{minipage}[h]{0.46\linewidth} \vskip 2mm
     System (\ref{eq:neural_model:ext}):
    $\thetavec=\col{-\frac{1}{\tau_m},\frac{A_f}{\tau_m},\frac{\sigma_s}{\tau_s}}$, $\lambda=\frac{\sigma_f}{A_f}$, $\bfA=\left(\begin{array}{cc} 0 & -\frac{1}{\tau_m}\\ 0 & -\frac{1}{\tau_s}\end{array}\right)$, $\Psivec(t,\lambda,y)=\left(\begin{array}{ccc} y & \tanh(\lambda y) & 0\\ 0&0& y\end{array}\right)$, $\bfg=0$
    \end{minipage}\\
  \hline
\end{tabular}
\end{table*}

For the time being we suppose that matrix $\bfA$ and vectors
$\bfB,\bfC$ in (\ref{eq:neural_model}) satisfy Assumption
\ref{assume:matrices} below.
\begin{assume}\label{assume:matrices} The triple  $\bfA,\bfB,\bfC$ is known, and there exist (and are known) a vector $\bfl$ and  matrices  $\bfP,\bfQ>0$ such that condition (\ref{eq:condition_MKY}) holds.
\end{assume}
 Note that Assumption \ref{assume:matrices} implies that the vector $\bfB=\mathrm{col}(1,b_1,\dots,b_{n-1})$ in (\ref{eq:neural_model}) is such that the polynomial $s^{n-1}+b_1s^{n-2}+\cdots+b_{n-1}$ is Hurwitz. At first, Assumption \ref{assume:matrices} may seem restrictive. In Section
\ref{sec:discussion} we lift this restriction by showing that the results presented for (\ref{eq:neural_model}) can be generalized to systems
\begin{equation}\label{eq:neural_model:ext}
\begin{array}{l}
\dot{\bfx}=\bfA \bfx + \Psivec(t,\lambdavec,y)\thetavec+\bfg(t,\lambdavec,y,u(t))+\xivec(t),\\
y=\bfC^{T}\bfx, \ \bfC=\col{1,0,\dots,0},
\end{array}
\end{equation}
in which the matrix $\bfA\in\Real^{n\times n}$ may be unknown but it is known that the pair $\bfA, \bfC$  is observable, the function
$\Psivec:\Real\times\Real^p\times\Real\rightarrow\Real^{n\times
m}$, $\Psivec\in\mathcal{C}^1$, is  Lipschitz in
$\lambdavec$,  and $\bfg(\cdot,\lambdavec,y(\cdot),u(\cdot))$, $\dot{\bfg}(\cdot,\lambdavec,y(\cdot),u(\cdot))$, ${\Psivec}(\cdot,\lambdavec,y(\cdot))$, $\dot{\Psivec}(\cdot,\lambdavec,y(\cdot))$ are bounded for all $\lambdavec\in\Omega_{\lambda}$ on $[t_0,\infty)$.

With regards to the functions $\varphivec$ and $\bfg$ in (\ref{eq:neural_model}) the following additional technical assumptions are made:

\begin{assume}\label{assume:varphi}
The functions $\varphivec(\cdot,\lambdavec,\cdot)$,
$\bfg(\cdot,\lambdavec,\cdot,\cdot)$ in (\ref{eq:neural_model})
are bounded and differentiable in
$\Real_{\geq t_0}\times\mathcal{D}_y$ and
$\Real_{\geq t_0}\times\mathcal{D}_y\times\mathcal{D}_u$ respectively,
and Lipschitz in $\lambdavec$. That is, there exist
 $D_{\varphi}, D_{g}, B_\varphi, B_g\in\Real_{\geq 0}$ such
that for all $t\in\Real_{\geq t_0}$, $y\in\mathcal{D}_y$, $u\in\mathcal{D}_u$,
$\lambdavec',\lambdavec''\in\Omega_{\lambda}$
\beq\label{eq:varphi_Lipschitz}
\begin{split}
    &\|\varphivec(t,\lambdavec',y)-\varphivec(t,\lambdavec'',y)\|\leq D_{\varphi} \|\lambdavec'-\lambdavec''\|,\\
    &\|\bfg(t,\lambdavec',y,u)-\bfg(t,\lambdavec'',y,u)\|\leq D_{g} \|\lambdavec'-\lambdavec''\|,
\end{split}
\eeq \beq\label{eq:boundedness_varphi_g}
  \|\varphivec(t,\lambdavec,y)\|\leq B_\varphi, \ \|\bfg(t,\lambdavec,y,u)\|\leq B_g.
\eeq Furthermore, there exist $M_\varphi,  M_g\in\Real_{\geq
0}$ such that
\begin{eqnarray}\label{eq:varphi_g_derivatives}
&&\begin{array}{l}\left|\frac{\pd \varphivec(t,\lambdavec,y)}{\pd y}\dot{y}+\frac{\pd \varphivec(t,\lambdavec,y)}{\pd t}\right|\leq M_\varphi\end{array}, \\
&&\begin{array}{l}\left|\frac{\pd \bfg(t,\lambdavec,y,u)}{\pd y}\dot{y}+\frac{\pd \bfg(t,\lambdavec,y,u)}{\pd
u}\dot{u}+\frac{\pd \bfg(t,\lambdavec,y,u)}{\pd t}\right|\leq M_g\end{array}\nonumber
\end{eqnarray}
for all $\lambdavec\in\Omega_\lambda$, $t\geq t_0$ along the solutions of
(\ref{eq:neural_model}).
\end{assume}
Conditions (\ref{eq:varphi_Lipschitz}), (\ref{eq:boundedness_varphi_g})
often hold naturally in the context of modeling and identification;
they may, however, impose limitations in the framework of
controller design. As for condition
(\ref{eq:varphi_g_derivatives}), the first inequality is a
version of (\ref{eq:derivative_bound}) that is essential for
uniform exponential convergence of solutions to the origin of
(\ref{eq:error_dynamics}) \cite{Automatica:Loria:2004}. The second
inequality in (\ref{eq:varphi_g_derivatives}) is a technical
condition.
Although this latter condition may
look somewhat restrictive, it may be relaxed if
$\bfg(t,\lambdavec,y,u(t))$ is expressed as
$\bfg(t,\lambdavec,y,u(t))=\bfg_1(t,y,u(t))+\bfg_2(t,\lambdavec,y,u(t))$.
In this case we would require that
(\ref{eq:boundedness_varphi_g}),
(\ref{eq:varphi_g_derivatives}) hold for $\bfg_2$.

A non-exhaustive list of systems that are
relevant in engineering applications and are governed by
(\ref{eq:neural_model}) or (\ref{eq:neural_model:ext}) includes
bio-/ chemical reactors \cite{Boskovic_1995,Comp:Poyton:2006},
nonlinear saturated magnetic circuits \cite{Moreau:2006},
magnetic bearings \cite{ACC:2000:Lin}, tire-road interaction,
and dynamics of live cells \cite{Rowat}. A few examples from this
list are provided in Table \ref{tab:examples}. The first model,
if described by (\ref{eq:neural_model}), trivially satisfies
Assumption \ref{assume:matrices}; it also satisfies Assumption
\ref{assume:varphi} if $T$ is bounded, differentiable,
separated away from zero, and  $y$, $\dot{T}$ are bounded.
 In the second model the pair $\bfA$, $\bfC$ is
observable, and $\varphivec$, $\dot{\varphivec}$ are bounded if
$y$, magnetic fluxes, expressed by $q_1$,$q_2$, and $\dot{q}_1$, $\dot{q}_2$ are bounded. This is
achievable via external controls \cite{ACC:2000:Lin}, at least
for small parametric mismatches  and $d$.  In the third model
the pair $\bfA$, $\bfC$ is observable, and boundedness of $y$,
$\dot{y}$, ${\Psivec}$, $\dot{\Psivec}$ is consistent with the
physics of the system.

\subsection{Problem formulation}

 Before we proceed with a formal problem statement, several points related to parametrization of (\ref{eq:neural_model}), (\ref{eq:neural_model:ext}) need to be discussed.
First, note that different definitions of systems (\ref{eq:neural_model}), (\ref{eq:neural_model:ext}) may correspond the same physical model. For example,  functions $\varphivec$ and $g$ and parameters $\thetavec, \lambdavec$ for the first model in Table \ref{tab:examples} may also be defined as:
\begin{eqnarray}\label{eq:param_ex_2}
&&\begin{array}{l}\varphivec(t,\lambdavec,y)=\col{1,y,e^{\lambda_1 \frac{1}{T(t)}+\lambda_2}y}, \ g=0, \end{array}
\end{eqnarray}
$\thetavec=\col{\frac{F}{V}x_0, -\frac{F}{V},-k_{\mathrm{ref}}}$, $\lambdavec=\col{-\frac{E}{R},\frac{E}{R T_{\mathrm{ref}}}}$, or
\begin{equation}\label{eq:param_ex_3}
\varphivec(t,\lambdavec,y)=\col{1,y}, \ g(t,\lambdavec,y,u(t))=-e^{\lambda_1 \frac{1}{T(t)}+\lambda_2}y,
\end{equation}
with $\thetavec=\col{\frac{F}{V}x_0, -\frac{F}{V}}$,  $\lambdavec=\col{-\frac{\ln(k_{\mathrm{ref}})E}{R},\frac{\ln(k_{\mathrm{ref}})E}{R T_{\mathrm{ref}}}}$. It is clear that if parametrization (\ref{eq:param_ex_2}) is chosen then identical outputs $y(t)$ will be observed for
 infinitely many combinations of parameters $\theta_3,\lambda_2$. Models of this type are referred to as {\it unidentifiable} \cite{Distefano:1980} (see also \cite{Chapell:1996},
\cite{Automatica:Denis-Vidal:2004}).
Dealing with unidentifiable models imposes technical difficulties. We will therefore assume that parametrizations which are obviously unidentifiable are avoided, if possible.
 As for the remaining alternative parametrizations, we assume that preference is given to those in which the dimension of $\lambdavec$ is minimal, i.e.  the parametrization in the first row in Table \ref{tab:examples} is preferable to (\ref{eq:param_ex_3}).

 Second, as far as identifiability is concerned, inferring true values of $\thetavec$, $\lambdavec$ from  output observations, $y(t)$, is not always possible, even if the system is linearly parameterized and no unmodeled dynamics are present. Consider
\begin{equation}\label{eq:example:non_identifiable}
\begin{split}
\dot{\bfx}&=\bfA \bfx+\left(\begin{array}{c}
1\\
1
\end{array}\right)\theta + \left(\begin{array}{c}
0\\
1
\end{array}\right)\lambda, \ \bfA=\left(\begin{array}{cc}-a_1 & 1\\ -a_2 & 0\end{array}\right)\\
y&=x_1, \ a_1,a_2\in\Real_{>0}.
\end{split}
\end{equation}
Let $\bfx(t,\theta,\lambda,\bfx_0)$ and
$\bfx(t,\theta',\lambda',\bfx_0')$ be two solutions of
(\ref{eq:example:non_identifiable}) corresponding to different
parameter values and initial conditions, and let
$\bfe=\mathrm{col}(e_1,e_2)=\bfx(t,\theta,\lambda,\bfx_0)-\bfx(t,\theta',\lambda',\bfx_0')$.
Picking $e_2(t_0)=-\theta+\theta'$, $e_1(t_0)=0$ ensures that
$e_1(t)=0$ for all $t\geq t_0$ if
$\theta-\theta'+\lambda-\lambda'=0$. Another, albeit nonlinearly parameterized, example is
\begin{equation}\label{eq:example:non_identifiable:2}
\dot{x}=-x + \theta + [\sin^2(\lambda+t)+x^2+1]^{-1}, \
 y=x.
\end{equation}
In this case
$x(t,\theta,\lambda,x_0)=x(t,\theta,\lambda',x_0)$ for all $t\geq t_0$ if $\lambda'=\lambda + k \pi$, $k\in\Numbers$, $\theta'=\theta$.

 In order to account for possible non-unique parametrization, for each pair $\thetavec,\lambdavec$ we introduce two sets: $\mathcal{E}_0(\lambdavec,\thetavec)$ and $\mathcal{E}(\lambdavec,\thetavec)$. The set $\mathcal{E}_0(\lambdavec,\thetavec)$:
\begin{eqnarray}\label{eq:goal_lambda_0}
&&\begin{array}{ll}
\mathcal{E}_0(\lambdavec,\thetavec)=&\{(\lambdavec',\thetavec'), \ \lambdavec'\in\Real^p, \thetavec'\in\Real^m| \\
 &\etavec_0(t,\lambdavec,\thetavec,\lambdavec',\thetavec')=0, \ \forall \ t\geq t_0\},
\end{array}\\
& &\etavec_0(t,\lambdavec,\thetavec,\lambdavec',\thetavec')= \nonumber \\
& & \ \ \ \ \ \ \  \ \ \ \ \bfB(\varphivec^{T}(t,\lambdavec,y(t)) \thetavec- \varphivec^{T}(t,\lambdavec',y(t))\thetavec')+\nonumber\\
&& \ \ \ \ \ \ \ \ \ \ \ \bfg(t,\lambdavec',y(t),u(t))-\bfg(t,\lambdavec,y(t),u(t))\nonumber
\end{eqnarray}
contains all  parametrizations of (\ref{eq:neural_model}) that are indistinguishable from observations of $\bfx(t)$ for the given $y(\cdot)$, $u(\cdot)$,  $\thetavec$, $\lambdavec$ at $\xivec(t)\equiv 0$. That is, if $\bfx(t,\thetavec,\lambdavec,\bfx_0)-\bfx(t,\thetavec',\lambdavec',\bfx_0)=0$ for all $t\geq t_0$ then $(\lambdavec',\thetavec')\in\mathcal{E}_0(\lambdavec,\thetavec)$. For system (\ref{eq:example:non_identifiable}) the set
$\mathcal{E}_0(\lambda,\theta)$ contains just one element, $(\lambda,\theta)$. The set $\mathcal{E}_0(\lambda,\theta)$, however, is not finite for system (\ref{eq:example:non_identifiable:2}) and for parametrization (\ref{eq:param_ex_2}) of the first model in Table \ref{tab:examples}.

The second set, $\mathcal{E}(\lambdavec,\thetavec)$, is defined as:
\begin{eqnarray}\label{eq:goal_lambda}
&&\mathcal{E}(\lambdavec,\thetavec)=\{(\lambdavec',\thetavec'),\lambdavec'\in\Real^p,\thetavec'\in\Real^m \ | \ \exists \ \bfp(\thetavec,\lambdavec,\thetavec',\lambdavec')\nonumber\\& &\in\Real^{n-1}:
\tilde{\bfC}^{T}e^{\Lambda(t-t_0)}\bfp +
  \etavec(t,\thetavec,\lambdavec,\thetavec',\lambdavec')=0 \ \ \forall \ t\geq t_0 \},\nonumber \\
\end{eqnarray}
where $\tilde{\bfC}\in\Real^{n-1}$, $\tilde{\bfC}=\mathrm{col}(1,0,\dots,0)$,
\begin{eqnarray}\label{eq:unidentifable_set_1}
& & \etavec(t,\lambdavec,\thetavec, \lambdavec',\thetavec') = \varphivec(t,\lambdavec,y(t))^{T}\thetavec -\varphivec(t,\lambdavec',y(t))^{T}\thetavec' + \nonumber\\
& & g_1(t,\lambdavec,y(t),u(t))-g_1(t,\lambdavec',y(t),u(t))+q(t,\lambdavec,\lambdavec')
\end{eqnarray}
and
$q(t,\lambdavec,\lambdavec')=\tilde{\bfC}^{T}\bfz(t,\lambdavec,\lambdavec')$:
\begin{eqnarray}\label{eq:unidentifable_set_2}
\dot\bfz&=& \Lambda \bfz +  \bfG (\bfg(t,\lambdavec,y(t),u(t))-\bfg(t,\lambdavec',y(t),u(t))), \nonumber \\
\Lambda&=&\left(\begin{array}{c|c}
\multirow{2}{*}{\ensuremath{-\bfb}} &   I_{n-2}\\ & 0 \end{array}\right), \ \bfG=\left(\begin{array}{cc} - {\bfb} & I_{n-1} \end{array}\right)\\
& & \bfz(t_0)=0, \ \bfb=\col{b_1,\dots,b_{n-1}}. \nonumber
\end{eqnarray}
If
$\xivec(t)\equiv0$ then the set $\mathcal{E}(\lambdavec,\thetavec)$
contains all  indistinguishable parametrizations of
(\ref{eq:neural_model})
for the given $y(\cdot)$, $u(\cdot)$, $\thetavec$, $\lambdavec$ (see Lemma
\ref{lem:observer_inferrence} in Section \ref{sec:proof}), and
$\mathcal{E}_0(\lambdavec,\thetavec)\subseteq\mathcal{E}(\lambdavec,\thetavec)$. Note that if $\mathrm{\dim}({\bfx})=2$ then $\Lambda=-b_1$,
$\bfG=(-b_1 \  1)$; if
$\mathrm{\dim}({\bfx})=1$ we set
$\mathcal{E}(\lambdavec,\thetavec)=\mathcal{E}_0(\lambdavec,\thetavec)$. For system
 (\ref{eq:example:non_identifiable}),
$\mathcal{E}(\lambda,\theta)=\{(\lambda',\theta'),
\lambda\in\Real, \ \theta\in\Real| \
\theta-\theta'+\lambda-\lambda'=0\}$.
Note that if $\bfg(t,\lambdavec,y(t),u(t))=\bfB g(t,\lambdavec,y(t),u(t))$, $\bfB=\mathrm{col}(1,b_1,\dots,b_{n-1})$, then $q(t,\lambdavec,\lambdavec')\equiv0$ for all $t\geq t_0$ in (\ref{eq:unidentifable_set_1}).  The introduction of sets $\mathcal{E}_0(\lambdavec,\thetavec)$, $\mathcal{E}(\lambdavec,\thetavec)$ does
not, of course, resolve identifiability issues. It helps, however, to specify
constraints on the nonlinearities in (\ref{eq:neural_model}) for which the parameter reconstruction, up to $\mathcal{E}_0$,
 $\mathcal{E}$, can be achieved.

Since this is desirable from an implementation point of view,
we will seek for a recursive procedure taking the values
$y(t)$, $u(t)$ as inputs and producing the estimates
$\hat{\bfx}(t)$, $\hat{\thetavec}(t)$, $\hat{\lambdavec}(t)$ of
$\bfx(t,\lambdavec,\thetavec,\bfx_0,[u])$, $\thetavec$, and
$\lambdavec$, respectively,  as outputs. Note that since we
allow for non-identifiable configurations, estimation of
parameters $\thetavec$, $\lambdavec$ is possible only up to
the set $\mathcal{E}(\lambdavec,\thetavec)$. We
will, therefore, be looking for an auxiliary system, i.e. an {\it adaptive
observer}: $\dot{\bfq}=\bff(t,y,u(t),\bfq)$,
$\bff:\Real\times\Real\times\Real\times\Real^{q}\rightarrow\Real^q$,
$\bfq(t_0)=\bfq_0$, $\bfq\in\Real^q$, and functions
$\bfh_{\theta}:\Real^q\rightarrow\Real^m$,
$\bfh_{\lambda}:\Real^q\rightarrow \Real^p$,
$\bfh_x:\Real^q\rightarrow\Real^n$ such that for the given $t_0$, appropriately chosen
$\bfq_0$,  known functions $r_1,r_2\in\mathcal{K}$, and all admissible $\lambdavec,\thetavec,\bfx_0$ the following requirements hold  for the
observer:
\begin{eqnarray}\label{eq:goal}
     & &\limsup_{t\rightarrow\infty}  \|\bfh_{x}(\bfq(t,\bfq_0))-\bfx(t,\lambdavec,\thetavec,\bfx_0)\|\leq r_1(\Delta_{\xi}) \\
    & &\limsup_{t\rightarrow\infty}   \
        \dist\left(\left(\begin{array}{c}\bfh_{\lambda}(\bfq(t,\bfq_0))\\
                                \bfh_{\theta}(\bfq(t,\bfq_0))
    \end{array}\right),\mathcal{E}(\lambdavec,\thetavec)\right)\leq r_2(\Delta_{\xi}),
    \nonumber
\end{eqnarray}
where $\Delta_\xi$ is defined in (\ref{eq:xi_bound}).


\section{Main result}\label{sec:results}

In this section we introduce an observer for \eqref{eq:neural_model}  and show that
%
asymptotic reconstruction of state and parameters of
\eqref{eq:neural_model}  is achievable (up to {the set $\mathcal{E}(\lambdavec,\thetavec)$}), subject to
some persistency of excitation conditions.

\subsection{Observer definition}\label{sec:definition_and_assumptions}

\begin{figure}[!t]
\centering
\includegraphics[width=0.7\linewidth]{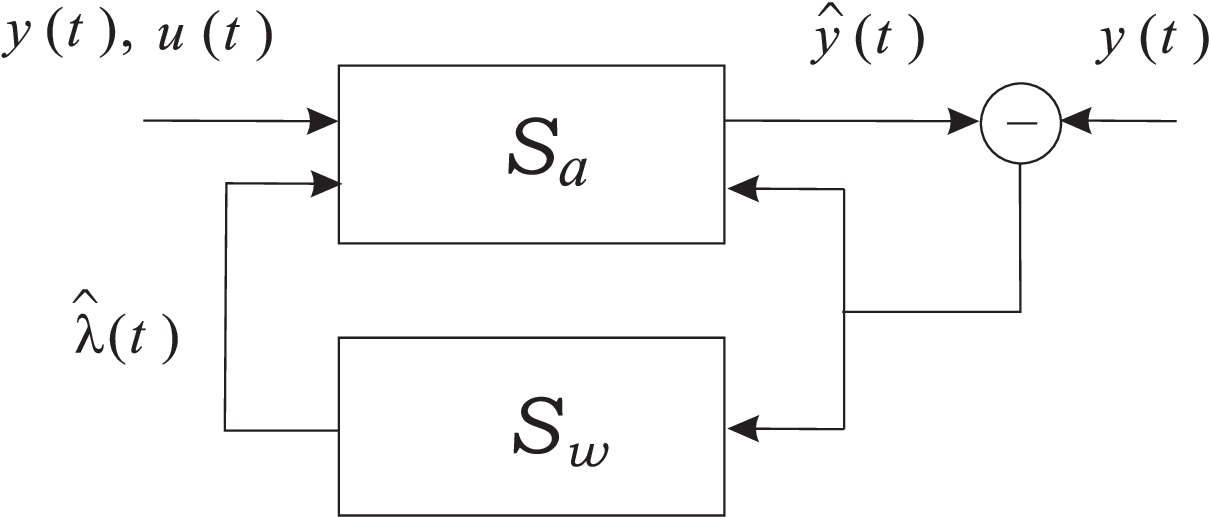}
\vspace{3mm}

\includegraphics[width=0.7\linewidth]{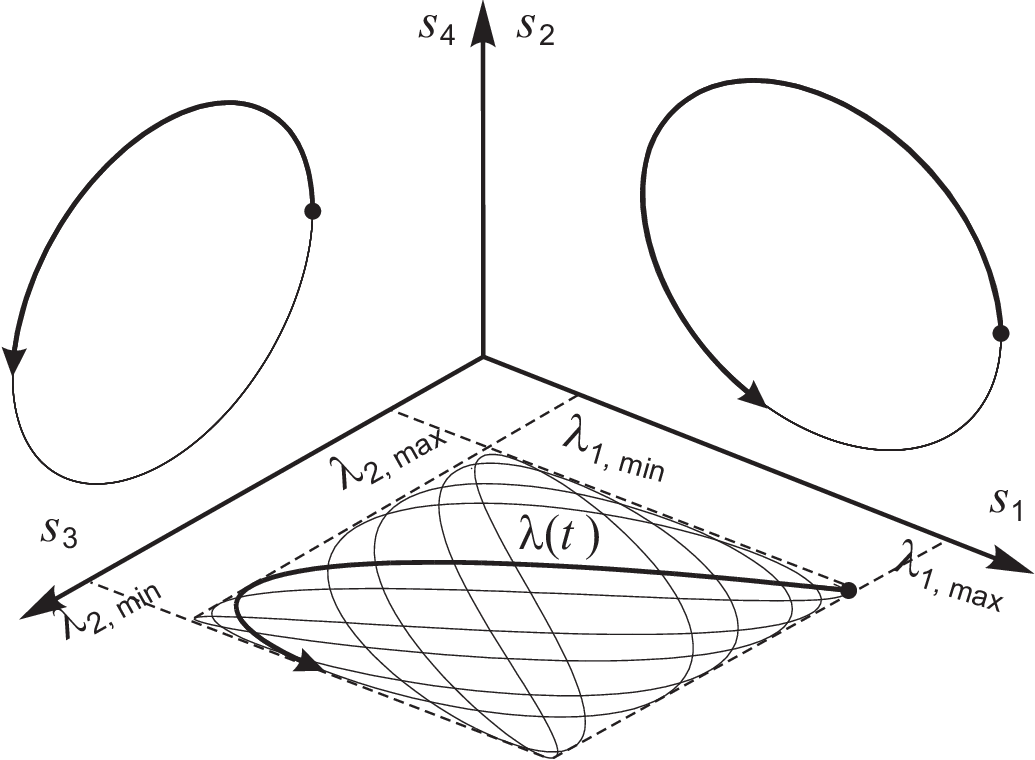}
    \caption{{\it Top panel}: general structure of the observer. {\it Bottom panel}: phase curves of system (\ref{eq:example_poisson}). }\label{fig:observer}
\end{figure}

Following the general ideas of
\cite{SIAM_non_uniform_attractivity} with regards to the
treatment of uncertain systems with general nonlinear
parametrization, we propose that an asymptotically converging
observer for (\ref{eq:neural_model}) consists of two coupled
subsystems, $\mathcal{S}_a$ and $\mathcal{S}_w$ (see Fig.
\ref{fig:observer}). The role of subsystem $\mathcal{S}_a$ is
to provide estimates of state and  parameters
$\thetavec$ of \eqref{eq:neural_model}, and the role of
subsystem $\mathcal{S}_w$ is to search the values of
 parameters $\lambdavec$.

The dynamics of subsystem $\mathcal{S}_a$ is defined as
follows:
\begin{eqnarray}\label{eq:linear_par_observer}
&&    \mathcal{S}_a:
        \left\{
            \begin{split}
             \dot{\hat{\bfx}} =& \bfA \hat{\bfx} +  \bfl (\bfC^{T}\hat{\bfx}-y(t)) + \bfB\varphivec^T(t,\hat{\lambdavec}(t),y(t)) \hat{\thetavec}\\ & + \bfg(t,\hat{\lambdavec}(t),y(t),u(t)) \\
                \dot{\hat{\thetavec}}=&-\gamma_{\theta} (\bfC^{T}\hat{\bfx}-y(t)) {\varphivec}(t,\hat{\lambdavec}(t),y(t)), \\
                \hat{y}=&\bfC^T \hat{\bfx}, \ \hat{\bfx}(t_0)\in\Real^n, \ \hat{\thetavec}(t_0)\in\Real^m,\ \gamma_{\theta}\in\Real_{>0}
            \end{split}
        \right.
\end{eqnarray}
The variable
$\hat{\thetavec}=\mathrm{col}(\hat{\theta}_1,\cdots,\hat{\theta}_m)$
in \eqref{eq:linear_par_observer} is an estimate of
$\thetavec$, and
$\hat{\lambdavec}=\mathrm{col}(\hat{\lambda}_1,\dots,\hat{\lambda}_p)$
is an estimate of $\lambdavec$. For the time being we suppose
that $\hat{\lambdavec}$  is a continuous function of $t$. The
matrix $\bfA$ and vectors $\bfB$, $\bfC$ in
(\ref{eq:linear_par_observer}) are identical to those in
(\ref{eq:neural_model}), and the vector $\bfl$ in
\eqref{eq:linear_par_observer} satisfies Assumption
\ref{assume:matrices}. If the values of $\lambdavec$
would be known then substitution $\hat{\lambdavec}=\lambdavec$
reduces system \eqref{eq:linear_par_observer} to
(\ref{eq:canonic_observer}), and conditions for asymptotic
reconstruction of state and parameter values of
(\ref{eq:neural_model}) follow from Theorem
\ref{thm:preliminary:exponential_skew_symmetric}. The values of $\lambdavec$, however, are unknown and therefore a
procedure for estimating the values of  $\lambdavec$ is needed.

Regarding the definition of $\mathcal{S}_w$, we propose that
the values of $\hat{\lambdavec}(t)$ result from an explorative
search in the domain $\Omega_{\lambda}$ of the admissible
values for $\lambdavec$.  The exploration can be realized by
movements along the solutions of a certain class of dynamical
systems. Let us, for instance, consider systems governed by the
following equations
\begin{equation}\label{eq:generator_poisson}
    \begin{split}
        \dot{\bfs}&=\bff(\bfs), \ \bfs(t_0)=\bfs_0, \
        \hat{\lambdavec}=\betavec(\bfs)
    \end{split}
\end{equation}
where $\bff:\Real^{n_p}\rightarrow\Real^{n_p}$,
$\betavec:\Real^{n_p}\rightarrow\Real^{p}$ are continuous, and
let $\Omega_s$ be the $\omega$-limit set\footnote{Recall that a
point $\bfz\in\Real^{n_p}$ is an $\omega$-limit point of
$\bfz_0\in\Real^{n_p}$ if there is a sequence $\{t_i\}$,
$i=1,2,\dots$, $\lim_{i\rightarrow\infty}t_i=\infty$, such that
$\lim_{i\rightarrow\infty} \bfs(t_i,\bfz_0)=\bfz$. The set of
all $\omega$-limit points of $\bfz_0$ is the $\omega$-limit set
of $\bfz_0$.} of $\bfs_0$. In addition, suppose that the
following properties hold:

\begin{enumerate}
\item[P1)] the functions $\bff$, $\betavec$ in
    (\ref{eq:generator_poisson}) are Lipschitz;
\item[P2)] the vector $\bfs_0$ is such that the solution
    $\bfs(t,\bfs_0)$ is bounded for all $t\geq t_0$;
\item[P3)] the image of $\Omega_s$ under $\betavec$
    contains $\Omega_\lambdavec$:  for every
    $\lambdavec\in\Omega_\lambda$ there is an
    $\bfs\in\Omega_s$ such that
    $\betavec(\bfs)=\lambdavec$.
\end{enumerate}

Properties P1 and P2 are technical requirements ensuring that
the derivative of $\hat\lambdavec$, as a function of $t$, is
bounded and has a bounded growth rate. Property P3, however, is
essential. It implies that the projection
$\betavec(\bfs(\cdot,\bfs_0))$ of the trajectory $\bfs(\cdot,\bfs_0)$
onto $\Omega_\lambda$ is dense and recurring in
$\Omega_\lambda$:
\begin{equation}\label{eq:dense_recurrent}
\begin{split}
     \forall \ \lambdavec\in\Omega_\lambda, \ &\forall \ \varepsilon\in\Real_{>0}, \ \forall \ t\geq t_0  \\
    & \exists \ t'>t: \ \|\lambdavec-\betavec(\bfs(t',\bfs_0))\|<\varepsilon.
\end{split}
\end{equation}
Indeed, let $\lambdavec'$ be an element from
$\Omega_{\lambda}$. Then according to P3 there is an
$\bfs'\in\Omega_s$: $\betavec(\bfs')=\lambdavec'$. Since
$\Omega_s$ is the $\omega$-limit set of $\bfs_0$, we can
conclude that there is a sequence $\{t_i\}$, $i=1,2,\dots$,
$\lim_{i\rightarrow\infty}t_i=\infty$, such that
$\lim_{i\rightarrow\infty}\bfs(t_i,\bfs_0)=\bfs'$. Finally,
using the continuity of $\betavec$  we arrive at
$\lim_{i\rightarrow\infty}\betavec(\bfs(t_i,\bfs_0))=\lambdavec'$.
In other words, for any $\lambdavec'\in\Omega_\lambda$ and
$\varepsilon>0$ there will exist a sequence of time instances
$t_i: \ \lim_{i\rightarrow\infty} t_i =\infty$ such that
$\|\hat{\lambdavec}(t_i)- \lambdavec'\|<\varepsilon$, and hence
(\ref{eq:dense_recurrent}) follows.  An example of a very
simple system possessing a solution $\bfs(t,\bfs_0)$ and an
output function $\betavec$  satisfying  properties P1--P3 for
$\Omega_\lambda=[-1,1]^2$ is
\begin{equation}\label{eq:example_poisson}
\begin{split}
&\dot{s}_1=-\sqrt{2} s_2, \ \dot{s}_2=\sqrt{2}s_1, \ \dot{s}_3=-s_4, \
\dot{s}_4=s_3,\\
&\betavec(\bfs(t,\bfs_0))=\mathrm{col}(s_1(t),s_3(t)), \bfs_0=\col{1,0,1,0}.
\end{split}
\end{equation}
Phase curves of (\ref{eq:example_poisson}) are shown in  Fig.
\ref{fig:observer}, bottom panel. Projections of the initial
segment of the trajectory are shown by thick lines. After
evolving beyond the initial segment, the values of
$\betavec(\bfs(t,\bfs_0))$ will densely fill the set
$[\lambda_{1,\min},\lambda_{1,\max}]\times[\lambda_{2,\min},\lambda_{2,\max}]=[-1,1]^2$,
cf. \cite{Nemytskii:1960}.

The problem with using (\ref{eq:generator_poisson}) directly as
an estimator for $\lambdavec$ is that exploration of the set
$\Omega_\lambda$ continues indefinitely. For the purposes
of observer design we need to ensure that exploration of
$\Omega_\lambda$ stops once a sufficiently small neighborhood
of the set $\mathcal{E}(\lambdavec,\thetavec)$ has been
reached. To enable this, the explorative subsystem must be
supplied with an error measure. A function of
$\|y(t)-\hat{y}(t)\|_{\varepsilon}$ is a possible
candidate for such a measure. Thus we replace the earlier
definition (\ref{eq:generator_poisson}) for
$\hat{\lambdavec}$ with the following:
\begin{equation}\label{eq:generator_poisson:2}
    \begin{split}
        \dot{\bfs}&=\gamma \sigma(\|y(t)-\hat{y}(t)\|_\varepsilon) \bff(\bfs), \ \varepsilon \in\Real_{\geq 0}, \ \gamma\in\Real_{>0},\\
        \hat{\lambdavec}&=\betavec(\bfs), \ \bfs(t_0)=\bfs_0,
    \end{split}
\end{equation}
where $\sigma:\Real_{\geq 0}\rightarrow\Real_{\geq 0}$ is a
bounded Lipschitz function:
\begin{equation}\label{eq:generator_poisson:3}
\begin{split}
\exists  \ D_\sigma, \ M_{\sigma}\in\Real_{>0}: &\\
  \sigma(\upsilon)\leq M_{\sigma},&  \ \sigma(\upsilon)\leq D_\sigma \upsilon \  \ \forall \  \upsilon\geq 0
\end{split}
\end{equation}
such that $\sigma(\upsilon)>0$ for $\upsilon>0$, and
$\sigma(0)=0$.
%
%

For the sake of simplicity and without loss of generality,
instead of dealing with general systems
(\ref{eq:generator_poisson:2}), we will focus on a specific
system of equations:
\begin{equation}\label{eq:nonlinear_par_observer}
\begin{split}
    \mathcal{S}_w: \
                &\left\{\begin{array}{ll}
                    \dot{{s}}_{2j-1}=&\gamma \sigma(\|y(t)-\hat{y}(t)\|_\varepsilon) \cdot\omega_j  \cdot ({s}_{2j-1} - {s}_{2j}  \\
                    & - {s}_{2j-1}({s}_{2j-1}^2+{s}_{2j}^2))
                    \\
                    \dot{{s}}_{2j}=&\gamma \sigma(\|y(t)-\hat{y}(t)\|_\varepsilon)\cdot\omega_j \cdot ({s}_{2j-1}+ {s}_{2j} \\
                    &- {s}_{2j}({s}_{2j-1}^2+{s}_{2j}^2) )
                    \\
                    \hat{\lambda}_j&=\beta_j(\bfs), \ j=\{1,\dots,p\},
                 \end{array}\right.\\
                 \beta_j(\bfs)&=\lambda_{j,\min}
                    +\frac{\lambda_{j,\max}-\lambda_{j,\min}}{2}({s}_{2j-1}+1)
            \end{split}
        \end{equation}
\begin{equation}\label{eq:initial_conditions}
     \bfs_0=\bfs(t_0): \ {s}_{2j-1}^2(t_0)+{s}_{2j}^2(t_0)=1,
\end{equation}
where $\sigma$ is a function satisfying
(\ref{eq:generator_poisson:3}).  Parameters
$\omega_j\in\Real_{>0}$ in (\ref{eq:nonlinear_par_observer})
are supposed to be {\it rationally-independent}:
\begin{equation}\label{eq:rational_independence}
   \begin{array}{l} \sum_{j=1}^p \omega_j k_j\neq 0, \ \forall \ k_j\in \Numbers. \end{array}
\end{equation}
Equations
(\ref{eq:nonlinear_par_observer})--(\ref{eq:rational_independence})
are straightforward generalizations from the example system of
which the phase curves are shown in Fig. \ref{fig:observer}. If
the term $\gamma\sigma(\|y(t)-\hat{y}(t)\|_\varepsilon)$ in the
right-hand side of (\ref{eq:nonlinear_par_observer}) is
substituted with $1$, these equations satisfy the requirements
P1--P3. Indeed, we can immediately see that in this case
 $(s_{2j-1}(t,\bfs_0),s_{2j}(t,\bfs_0))=(\cos(\omega_j(t-t_0)  +
a_j),\sin(\omega_j(t-t_0)+a_j)), \ a_j\in\Real$. Thus
properties P1, P2 hold. Trajectories $s_{1}(\cdot,\bfs_0)$,
$s_{3}(\cdot,\bfs_0)$, $\dots$ $s_{2p-1}(\cdot,\bfs_0)$ evolve on a
corresponding $p$-dimensional invariant torus. Since $\omega_j$
are rationally-independent these trajectories densely fill the
torus (cf. \cite{Arnold_78}, \cite{Nemytskii:1960},
\cite{Katok:1999}) or, alternatively, they densely fill the
hypercube $[-1,1]^p$. This implies that $\Omega_s$, the
$\omega$-limit set of $\bfs_0$, is
$\Omega_s=\{\col{s_1,s_2,\dots,s_{2p}}\in\Real^p| \
\col{s_1,s_3,\dots,s_{2p-1}}\in[-1,1]^p, \
s_{2j}=\pm\sqrt{1-s_{2j-1}^2}, \ j=1,\dots,p\}$. Noticing that
the image of $\Omega_s$ under transformation  $\betavec$
coincides with $\Omega_\lambda$ we conclude that P3 holds.

Concerning the structure of $\mathcal{S}_w$, no additional
model-dependent constraints are imposed on
(\ref{eq:nonlinear_par_observer}) (or, in general, on
(\ref{eq:generator_poisson:2})), apart from the general
requirements P1--P3. Model-specific nonlinearities are
accounted for in the ``converging'' part, $\mathcal{S}_a$, of
the observer producing the estimates for $\thetavec$ and
$\bfx$. The information about the values of $\lambdavec$ is
transferred to the exploratory part, $\mathcal{S}_w$, by means
of $\|y(t)-\hat{y}(t)\|_\varepsilon$. The latter variable modulates
the speed of exploration in $\Omega_\lambda$ along a  search
trajectory. The  search trajectory itself does not need to be
dependent on the properties of $\bfg$, $\varphivec$, and
neither is the structure of $\mathcal{S}_w$. This potential advantage of the approach, however, comes at a cost.
According to
(\ref{eq:dense_recurrent}), small neighborhoods of sets to
which the solutions of the combined system
(\ref{eq:linear_par_observer}),
(\ref{eq:nonlinear_par_observer}) converge are not necessarily
forward invariant. Hence these sets are not guaranteed
to be asymptotically stable. Nevertheless, albeit in a weaker sense \cite{Milnor_1985}, they are still attracting. We illustrate this point in Section \ref{sec:example}.

\subsection{Asymptotic properties of the observer}\label{sec:main_result_observer}

Let us now proceed with specifying those properties of
\eqref{eq:neural_model} that can be useful for state and
parameter reconstruction. Recall that
(\ref{eq:neural_model}) is a generalization of the standard
canonic observer form (\ref{eq:observer_canonical_form}).
According to Theorem
\ref{thm:preliminary:exponential_skew_symmetric}, one of the
conditions for (\ref{eq:canonic_observer}) to be an adaptive
observer for (\ref{eq:observer_canonical_form}) is persistency
of excitation of the restriction of $\phivec(\cdot,y(\cdot))$ on $\Real_{\geq t_0}$. It is
therefore natural to expect that some kind of persistency of
excitation conditions might be needed for reconstruction of
parameters $\thetavec$, $\lambdavec$ in (\ref{eq:neural_model})
too. Two versions of these conditions will be considered,
namely the notions of {\it uniform persistency of excitation}
\cite{Lorea_2002} and {\it nonlinear persistency of excitation}
\cite{Cao_2003}.
\begin{defn}\label{defn:uniform_pe}
       A function $\alphavec:\Real_{\geq t_0} \times \Omega_\lambda \ra \Real^p$ is  $\lambda$-Uniformly Persistently Exciting ($\lambda$-UPE with $T$, $\mu$), denoted by $\alphavec(t,\lambdavec)\in\lambda{\rm UPE}(T,\mu)$, if there exist $T,\mu\in\Real_{>0}$:
    \beq\label{eq:PE_linear_uniform}
        \begin{array}{l}\int_{t}^{t+T} \alphavec(\tau,\lambdavec)\alphavec^T(\tau,\lambdavec){\rm d}\tau \geq \mu
         \bfI_p, \ \forall \ t\geq t_0, \ \lambdavec\in
        \Omega_\lambda.\end{array}
    \eeq
\end{defn}
In contrast to the conventional definitions of persistency of
excitation (cf. Definition \ref{defn:pe}),  uniform persistency
of excitation requires existence of $\mu,T\in\Real_{>0}$ in
(\ref{eq:PE_linear_uniform}) that are independent on
$\lambdavec$ for all $\lambdavec\in\Omega_\lambda$. This is a
stronger restriction; we will, however, require that it holds for
$\varphivec(t,\lambdavec,y(t))$ (as a function of $t,\lambdavec$ on $\Real_{\geq t_0}\times\Omega_\lambda$) in
\eqref{eq:neural_model}.

Since parametrization of (\ref{eq:neural_model}) is allowed to
be nonlinear, it is natural to expect that reconstruction of
model parameters might require a nonlinear version of standard
persistency of excitation. Here we employ the following
generalization of the standard notion (cf. \cite{Cao_2003}):
\begin{defn}\label{defn:nonlinear_pe} Let $\mathcal{E}$ be
a set-valued map defined on $\mathcal{D}\subset\Real^d$  and
associating a subset of $\mathcal{D}$ to every
$\bfp\in\mathcal{D}$.
   A function $\alphavec:\Real_{\geq t_0} \times \mathcal{D}\times\mathcal{D} \ra \Real^k$ is
    said to be weakly Nonlinearly Persistently Exciting in $\bfp$ wrt $\mathcal{E}$ (wNPE with $L,\beta,\mathcal{E}$ ),
    denoted by $\alphavec(t,\bfp,\bfp')\in{\rm wNPE}(L,\beta,\mathcal{E})$, if there exist $L\in\Real_{>0}$, $t_1\geq t_0$, and $\beta\in\mathcal{K}_{\infty}$:
    \beq\label{eq:nonlinear_pe}
    \begin{split}
     & \forall \ t\geq t_1, \ \bfp,\bfp'\in\mathcal{D}  \ \ \exists \  t'\in[t,t+L]:  \\
     & \ \ \ \ \ \ \ \ \norm{\alphavec(t',\bfp,\bfp')}
     \geq \beta \left({\rm dist}(\Ecal(\bfp),\bfp^\prime)\right).
    \end{split}
    \eeq
\end{defn}
If the set  $\mathcal{E}(\bfp)$ contains just one element,
$\bfp$, then  the inequality in (\ref{eq:nonlinear_pe}) reduces
to $\norm{\alphavec(t',\bfp,\bfp')} \geq \beta
(\|\bfp-\bfp^\prime\|)$.   Taking the above notions into account we formulate the main technical assumption on the nonlinearities in (\ref{eq:neural_model}):

\begin{assume}\label{assume:nonlinearities} The functions $\varphivec$, $\bfg$ in the right-hand side of (\ref{eq:neural_model}) are such that
\begin{enumerate}
\item[ A1) ] the restriction  of the function $\alphavec_1:\Real\times\Real^p\rightarrow\Real^m$,
    $\alphavec_1(t,\lambdavec)=\varphivec(t,\lambdavec,y(t))$ on
    $\Real_{\geq t_0}\times\Omega_{\lambda}$ is $\lambda$-UPE with $T,\mu$;
\item[ A2) ] the restriction of $\alphavec_2:\Real\times\Real^{p+m}\times\Real^{p+m}\rightarrow\Real$, $\alphavec_2(t,(\lambdavec,\thetavec),(\lambdavec',\thetavec'))=\etavec(t,\lambdavec,\thetavec,\lambdavec',\thetavec')$,
    where $\etavec(\cdot)$ is defined in (\ref{eq:unidentifable_set_1}), on $\Real_{\geq t_0}\times\Real^{p+m}\times\Real^{p+m}$ is weakly nonlinearly persistently exciting in $(\lambdavec,\thetavec)$ wrt to the map $\mathcal{E}(\lambdavec,\thetavec)$ determined by (\ref{eq:goal_lambda}).
\end{enumerate}
\end{assume}
\begin{rem}\label{rem:relax_pe}\normalfont
Checking that condition A1 holds is straightforward if e.g.  $\varphivec(t,\lambdavec,y(t))$ is periodic in $t$. Regarding condition A2, we note that, according to (\ref{eq:unidentifable_set_1}),  $\etavec(t,\lambdavec,\thetavec,\lambdavec',\thetavec')$ can be expressed as
$\etavec(t,\lambdavec,\thetavec,\lambdavec',\thetavec')= r(t,\lambdavec,\thetavec)- r(t,\lambdavec',\thetavec)+ \varphivec(t,\lambdavec',y(t))^{T}(\thetavec-\thetavec')$,
where
\[
\begin{array}{l}
r(t,\lambdavec,\thetavec)=\varphivec(t,\lambdavec,y(t))^{T}\thetavec+g_1(t,\lambdavec,y(t),u(t))\\
+\tilde{\bfC}^{T}\int_{t_0}^{t}e^{\Lambda(t-\tau)}\bfG \bfg(\tau,\lambdavec,y(\tau),u(\tau))d\tau.
\end{array}
\]
 If $\varphivec,\bfg$ are differentiable then
$r(t,\lambdavec,\thetavec)- r(t,\lambdavec',\thetavec)=
\bfR(t,\lambdavec,\lambdavec',\thetavec)$
$(\lambdavec-\lambdavec')$, where $\bfR(t,$
$\lambdavec,\lambdavec',$ $\thetavec)$ $=$ $\int_{0}^1
\frac{\pd }{\pd \bfs} r(t$
$,\bfs(\xi,\lambdavec,\lambdavec'),\thetavec) d\xi$,
$\bfs(\xi,\lambdavec,\lambdavec')=\lambdavec \xi +
(1-\xi)\lambdavec'$. Hence
\begin{equation}\label{eq:nonlinear_pe_test}
\begin{array}{l}
\etavec(t,\lambdavec,\thetavec,\lambdavec',\thetavec')=\\
\ \ \ (\varphivec^T (t,\lambdavec',y(t)),\bfR(t,\lambdavec,\lambdavec',\thetavec))(\thetavec-\thetavec',\lambdavec-\lambdavec')
\end{array}
\end{equation}
It is therefore clear that if  $(\varphivec^T(t,\lambdavec',y(t)),\bfR(t,\lambdavec,\lambdavec',\thetavec))$, $t\geq t_0$, $\lambdavec,\lambdavec'\in\Omega_\lambda$, $\thetavec\in\Omega_\theta$, is $(\lambda,\lambda',\theta)$-uniformly persistently exciting,  then
the system is uniquely identifiable, it satisfies condition A2, and $\mathcal{E}_0(\lambdavec,\thetavec)$ coincides with $\mathcal{E}(\lambdavec,\thetavec)$.
\end{rem}
We are now ready to state the main
result:
\begin{thm}\label{theorem:neural_identification}
    Consider  system \eqref{eq:neural_model} together with the observer defined by \eqref{eq:linear_par_observer}, \eqref{eq:nonlinear_par_observer}--\eqref{eq:rational_independence}. Suppose that Assumptions \ref{assume:matrices},
    \ref{assume:varphi}, and \ref{assume:nonlinearities}  hold. Then there exist a constant $\bar \gamma \in \Real_{>0}$ and
functions $r_1,r_2\in\mathcal{K}$ such that if $\gamma,
\varepsilon$ are the corresponding parameters of
(\ref{eq:nonlinear_par_observer}), and $\gamma \in (0, \bar
\gamma)$, $\varepsilon > r_1(\Delta_{\xi})$,  then
\begin{equation}\label{eq:thm2}
        \begin{array}{c}\limsup_{t\ra\infty}
        \dist\left(\left(\begin{array}{c}\hat\lambdavec(t)\\
                                \hat\thetavec(t)
    \end{array}\right),\mathcal{E}(\lambdavec,\thetavec)\right)\leq r_2(\varepsilon).
\end{array}
\end{equation}
If, in addition,  $\mathcal{E}(\lambdavec,\thetavec)$ coincides
with $\mathcal{E}_0(\lambdavec,\thetavec)$ then there is an
$r_3\in\mathcal{K}$:
\begin{align}
        &\limsup_{t\ra\infty}\norm{\hat \bfx(t)-\bfx(t)} \leq r_3(\varepsilon). \label{eq:thm1}
\end{align}
\end{thm}

The proof of Theorem \ref{theorem:neural_identification} is
presented in the next section.

Let us briefly comment on the assumptions made in the
theorem statement. Assumptions \ref{assume:matrices}, \ref{assume:varphi} are standard; condition  A1 in Assumption \ref{assume:nonlinearities} is the
conventional requirement ensuring exponential convergence of
$\hat \bfx(t)$, $\hat \thetavec(t)$ to $\bfx(t)$ and
$\thetavec$ provided that the value of $\lambdavec$ and hence
the values of $\varphivec(t,\lambdavec,y(t))$ are known (cf.
Theorem \ref{thm:preliminary:exponential_skew_symmetric});
 A2 ensures that the distance from
$(\hat{\lambdavec},\hat{\thetavec})$ to the set
$\mathcal{E}(\lambdavec,\thetavec)$ is inferable from the values
of $y(t)-\hat{y}(t)$ over $[t_0,\infty)$ (cf. Lemma
\ref{lem:observer_inferrence}). Note that nonlinear dependence
of $\etavec$ on $\lambdavec,\lambdavec'$ can impose certain
technical and computational difficulties whilst checking that
this condition holds. Finally, observe that the state estimation in
the proposed scheme requires that
$\mathcal{E}_0(\lambdavec,\thetavec)=\mathcal{E}(\lambdavec,\thetavec)$.
System
(\ref{eq:example:non_identifiable}) illustrates that violation of
this assumption may prevent the reconstruction of the state from
observed output data.

The value of $\bar{\gamma}$ and the functions $r_1,r_2,r_3$ could
in principle be given explicitly. However, due to dependence of
these functions on $\bfA$, $\bfB$, $\bfC$, parameters
$D_\varphi$, $D_g$, $M_\varphi$, $M_g$  and $T$, $L$, $\mu$,
$\beta$ specified in Assumptions \ref{assume:matrices},
\ref{assume:varphi} and Definitions \ref{defn:uniform_pe},
\ref{defn:nonlinear_pe},  explicit expressions for
$\bar{\gamma}$ and $r_1,r_2,r_3$  are too lengthy and
thus are removed from the theorem's statement. They are,
nevertheless, provided in the proof (see e.g.
(\ref{eq:epsilon_gamma_choice}), (\ref{eq:r_1_2_estimate})). A procedure for
finding the values of $\bar{\gamma}$ and $\varepsilon$ is
discussed in Section \ref{sec:example}.

The value of $\varepsilon$, viz. the accuracy of the
estimation, is determined by the bound $\Delta_\xi$ on the
amplitude of perturbation $\xivec(t)$. This dependency is
established through the function $r_1$ determining
a lower bound for parameter $\varepsilon$ in
\eqref{eq:nonlinear_par_observer}. If no perturbation
$\xivec(t)$ is present in the right-hand side of
(\ref{eq:neural_model}) then the value of $\varepsilon$ can be
chosen arbitrarily small. Note that the convergence itself is
asymptotic and not necessarily exponential. This is the price
for the presence of unknown parameters $\lambdavec$ in
(\ref{eq:neural_model}).

\begin{rem}\label{rem:jumps}\normalfont
The estimate $\hat{\lambdavec}(t)$ is guaranteed to converge to a single element of $\Omega_\lambda$ (see (\ref{eq:convervence_lambda_h}));
estimates $\hat{\thetavec}(t)$ may oscillate due the influence of $\xivec(t)$. These oscillations are bounded, and will eventually be confined to the $2 r_2(\varepsilon)$-neighborhood of $\mathcal{E}(\lambdavec,\thetavec)$. Hence, for $t_1$ sufficiently large and all $t\geq t_1 \geq t_0$, the $2 r_2(\varepsilon)$-neighborhood of $(\hat{\lambdavec}(t),\hat{\thetavec}(t))$ will always contain an element of $\mathcal{E}(\lambdavec,\thetavec)$.  This element may not necessarily be from $\Omega_\lambda\times\Omega_\theta$.  If the elements of $\mathcal{E}(\lambdavec,\thetavec)$ are separated by distances exceeding $3 r_2(\varepsilon)$ then the estimates are guaranteed to converge to  the $r_2(\varepsilon)$-vicinity of just one element. This point in $\mathcal{E}(\lambdavec,\thetavec)$ will depend on $\xivec$, $\bfx_0$, and on the initial state of the observer.
\end{rem}

\begin{rem}\label{rem:relax}\normalfont The function $\beta$ in Definition \ref{defn:nonlinear_pe} can be allowed to depend
on $(\lambdavec,\thetavec)$. In view of Remark \ref{rem:relax_pe}, this relaxes the requirement that
 $(\varphivec^{T}(t,\lambdavec',y(t)),\bfR(t,\lambdavec,\lambdavec',\thetavec))$ in (\ref{eq:nonlinear_pe_test}) (as a function of $t$, $(\lambdavec,\lambdavec',\thetavec)$ for $t\geq t_0$)
is $(\lambda,\lambda',\theta)$-UPE to the condition that $(\varphivec^{T}(t,\lambdavec',y(t))$ $,$ $\bfR(t,\lambdavec,\lambdavec',\thetavec))$ is  $\lambda'$-UPE.  Note that this will make $r_2$ in (\ref{eq:thm2}) dependent on $(\lambdavec, \thetavec)$. Finally, note that A2 need not hold for all $(\lambdavec,\thetavec)\in\Real^{p+m}$ and can be restricted to the union of $\Omega_\lambda\times\Omega_\theta$ and the domain to which $(\hat{\lambdavec}(t),\hat{\thetavec}(t))$ belong for $t\geq t_0$.
\end{rem}


\section{Proof of Theorem \ref{theorem:neural_identification}}\label{sec:proof}

According to Assumption \ref{assume:varphi} and
(\ref{eq:xi_bound}),  functions $\varphivec$, $\bfg$ and
$\xivec$ are continuous and are bounded  in $\Real_{\geq t_0}\times\Omega_\lambda\times \mathcal{D}_y$, $\Real_{\geq t_0}\times\Omega_\lambda\times \mathcal{D}_y\times\mathcal{D}_u$ and $\Real_{\geq t_0}$ respectively. Therefore solutions of the
combined system, (\ref{eq:neural_model}),
(\ref{eq:linear_par_observer}),
(\ref{eq:nonlinear_par_observer})--(\ref{eq:rational_independence})
exist and are defined for all $t\geq t_0$.  Let us denote
$\bfe=\mathrm{col}(\bfe_1,\bfe_2)$,  $\bfe_1:=\hat \bfx -
\bfx$, $\bfe_2:=\hat \thetavec - \thetavec$, $\alphavec(t,\hat
\lambdavec)=\varphivec(t,\hat \lambdavec,y(t))$. Then according
to \eqref{eq:neural_model} and \eqref{eq:linear_par_observer}
the following holds along the solutions of
(\ref{eq:neural_model}), (\ref{eq:linear_par_observer}),
(\ref{eq:nonlinear_par_observer})--(\ref{eq:rational_independence})
\begin{eqnarray}\label{eq:error_sys}
       \bpm \dot \bfe_1 \\ \dot \bfe_2 \epm
        &=&
        \bpm
            \bfA+\bfl \bfC^T   &   \bfB \alphavec^T(t,\hat \lambdavec(t)) \\
            -\gamma_{\theta} \alphavec(t,\hat \lambdavec(t))\bfC^T & 0
        \epm
        \bpm
            \bfe_1 \\ \bfe_2
        \epm        \nonumber \\
        & & +
        \bpm
            \bfv(t,\hat \lambdavec(t),\lambdavec, y(t),u(t)) \\ 0
        \epm
\end{eqnarray}
where the function $\bfv$:
\beq\label{eq:defn_v}
\begin{split}
    &\bfv(t,\hat\lambdavec, \lambdavec,y,u)  = \bfB(\varphivec^T(t,\hat \lambdavec,y)-\varphivec^T(t,\lambdavec,y))\thetavec \\
    & +  \bfg(t,\hat{\lambdavec},y,u)-\bfg(t,\lambdavec,y,u) - \xivec(t).
\end{split}
\eeq is continuous and bounded  for all $y\in\mathcal{D}_y$, $u\in\mathcal{D}_u$, $\lambdavec,\hat{\lambdavec}\in\Omega_\lambda$ and $t\geq t_0$.

The proof of the theorem  is split into three parts. In the {\it
first}  part we consider systems
 \beq
\label{eq:sys_lemma_exponential}
            \bpm \dot \bfe_1 \\ \dot \bfe_2 \epm
            =
            \bpm
                \bfA+\bfl \bfC^T   &   \bfB \alphavec^T(t,\hat \lambdavec(t)) \\
                -\gamma_{\theta} \alphavec(t,\hat \lambdavec(t))\bfC^T & 0
            \epm
            \bpm
                \bfe_1 \\ \bfe_2
            \epm,
        \eeq
where  $\gamma_{\theta}\in\Real_{>0}$ and
 $\hat \lambdavec:\Real_{\geq t_0}\rightarrow\Omega_{\lambda}$ is a continuous, differentiable and bounded function.
 Let $\Phi(t)$ be a fundamental solution matrix of (\ref{eq:sys_lemma_exponential}), and denote
$\Phi(t,t_0)=\Phi(t)\Phi(t_0)^{-1}$.
Since $\Phi(t_0,t_0)$ is the identity matrix we say that
$\Phi(t,t_0)$ is the normalized solution matrix of
(\ref{eq:sys_lemma_exponential}). We show that if  Assumptions \ref{assume:matrices}, \ref{assume:varphi} and condition A1 of Assumption \ref{assume:nonlinearities}
hold then there are positive numbers
$M_\lambda$, $\rho$, and $D_\rho$ such that the fundamental
(normalized) matrix of solutions, $\Phi(t,t_0)$, of
(\ref{eq:sys_lemma_exponential}) with $\hat \lambdavec: \
\|\dot{\hat \lambdavec}(t)\|\leq M_\lambda$ satisfies the
inequality $\|\Phi(t,t_0) \bfp \|\leq D_\rho
e^{-\rho(t-t_0)}\|\bfp\|, \ \bfp\in\Real^{n+m}$, $t\geq t_0$.

Using this result, in the {\it second part} of the proof we
demonstrate existence of $\bar \gamma$ and an $\varepsilon$, dependent on $\Delta_{\xi}$,
such that setting $\gamma\in(0,\bar \gamma]$
ensures that the
estimate $\hat{\lambdavec}(t)$ converges to a
$\lambdavec^{\ast}$ from $\Omega_\lambda$:
$\lim_{t\rightarrow\infty}\hat{\lambdavec}(t)=\lambdavec^\ast$.

{ In the {\it third  part} of the proof we show that   condition A2 of Assumption \ref{assume:nonlinearities} guarantees
that (\ref{eq:thm1}) holds and that, subject to the condition that $\mathcal{E}_0(\lambdavec,\thetavec)$ coincides with $\mathcal{E}(\lambdavec,\thetavec)$, property (\ref{eq:thm2}) holds too.}

{\it Part 1} of the proof is contained in the
following.

\begin{lem}\label{lem:UPE} Consider  system (\ref{eq:sys_lemma_exponential}), and suppose that
    \begin{enumerate}
        \item[C1) ] the  matrices $\bfA$, $\bfB$, and
            $\bfC$ satisfy Assumption \ref{assume:matrices}
        \item[C2) ] the restriction of the function $\alphavec$ in the
            right-hand side of
            (\ref{eq:sys_lemma_exponential}) on $\Real_{\geq t_0}\times\Omega_\lambda$ is
            $\lambda$-UPE with constants $T,\mu$ as in
            \eqref{eq:PE_linear_uniform}
        \item[C3) ] the function
            $\alphavec(t,\cdot)$ is Lipschitz in
            $\Omega_{\lambdavec}$ uniformly in $t$, $t\geq t_0$: there is a $D\in\Real_{>0}$ such that
            $\|\alphavec(t,\lambdavec)-\alphavec(t,\lambdavec')\|\leq
            D \|\lambdavec-\lambdavec'\| \ \forall  \
            \lambdavec,\lambdavec'\in\Omega_{\lambda}$, $t\geq t_0$
        \item[C4) ] the function $\alphavec$ and its
            partial derivatives wrt $t$, $\lambdavec$ are
            bounded; that is there is a constant $M$ such
            that
            $\max\left\{\left\|\alphavec(t,\lambdavec)\right\|,\left\|\frac{\pd
            \alphavec(t,\lambdavec)}{\pd t}\right\|,
            \left\|\frac{\pd \alphavec(t,\lambdavec)}{\pd
            \lambdavec}\right\| \right\}\leq M \ \forall \
            \lambdavec\in\Omega_{\lambda}, \ \forall \
            t\geq t_0$.
         \end{enumerate}

Then there exist $\rho, D_{\rho}\in\Real_{>0}$ such
that for all $\hat \lambdavec:\Real_{\geq t_0}\rightarrow\Omega_\lambda$, $\hat \lambdavec\in\mathcal{C}^1$:
\begin{equation}\label{eq:lambda_condition}
\|\dot{\hat \lambdavec} (t)\|\leq M_\lambda
\end{equation}
\begin{equation}\label{eq:gamma_boundary:1}
0 \leq M_{\lambda} \leq {\mu r}/({2 D M T^2}), \ r\in[0,1)
\end{equation}
the following holds
  \beq\label{eq:exponential_UPE}
        \norm{\Phivec(t_2,t_1)\bfp} \leq e^{-\rho(t_2-t_1)}D_\rho\norm{\bfp}, \ \bfp\in\Real^{n+m},
  \eeq
where $t_2\geq t_1\geq t_0$, and $\Phivec(\cdot,\cdot)$ is the
fundamental (normalized) solution matrix of
\eqref{eq:sys_lemma_exponential}.
\end{lem}
The proof of Lemma \ref{lem:UPE} and other auxiliary results are
provided in Appendix.

{\it Part 2}. Consider now the interconnection of
(\ref{eq:neural_model}) with the observer
(\ref{eq:linear_par_observer}),
(\ref{eq:nonlinear_par_observer})--(\ref{eq:rational_independence}).
The dynamics of the combined system in the coordinates $\bfe$,
$\hat{\lambdavec}$ is described by  (\ref{eq:error_sys}),
(\ref{eq:nonlinear_par_observer})--(\ref{eq:rational_independence}).
Recall that $\bfe(t)$, $\hat{\lambdavec}(t)$ are defined for
all  $t\geq t_0$. With respect to $\bfe(t)$, the following
holds
\begin{equation}\label{eq:error_dynamics:explicit}
\begin{array}{l}
\bfe(t)=\Phi(t,t_0)\bfe(t_0)\\
+\int_{t_0}^{t}\Phi(t,\tau) \bpm
            \bfv(\tau,\hat \lambdavec(\tau),\lambdavec, y(\tau),u(\tau)) \\ 0
        \epm d\tau.
\end{array}
\end{equation}
The variable $\hat \lambdavec$ in the combined system, as a
function of $t$, is bounded, continuous, and differentiable
with bounded derivative. Moreover, for any $\bfA,\bfB,\bfC$,
$\thetavec\in\Omega_{\theta}$, $\lambdavec\in\Omega_{\lambda}$
and for any choice of $\bfl$, $\gamma_\theta$ in
(\ref{eq:linear_par_observer}) we have that $|\dot {\hat
\lambda}_j(t)|\leq \gamma  M_\sigma \max_{i} |\omega_{i}|
\frac{\lambda_{i,\max}-\lambda_{i,\min}}{2}$, where $M_\sigma$,
$\omega_i$, and $\gamma$ are defined in
(\ref{eq:generator_poisson:3}),
(\ref{eq:rational_independence}), and
(\ref{eq:nonlinear_par_observer}). Thus   $\forall \ t\geq t_0$
we have:
\begin{equation}\label{eq:lambda_condition:1}
\begin{array}{l}
\hat{\lambdavec}(t)\in\Omega_{\lambda}, \ \|\dot{\hat \lambdavec}(t)\| \leq \gamma \sqrt{p} M_\sigma \max_{i} |\omega_{i}| \frac{\lambda_{i,\max}-\lambda_{i,\min}}{2},
\end{array}
\end{equation}
where $p$ is the dimension of $\lambdavec$. According to
 Assumption \ref{assume:varphi},  function $\alphavec$ is bounded, differentiable, and Lipschitz in the
second argument.  In particular, $\|\alphavec(t,\hat
\lambdavec)\|\leq B_\varphi, \ \left\|\frac{{\alphavec}(t,\hat
\lambdavec)}{\pd t}\right\|\leq M_\varphi, \
\left\|\frac{{\alphavec}(t,\hat \lambdavec)}{\pd \hat
\lambdavec}\right\|\leq D_\varphi$  $\forall$ $t\geq t_0$, $\hat{\lambdavec}\in\Omega_{\lambda}$. Hence there is an
$M=\max\{B_\varphi,M_\varphi,D_\varphi\}$ such that condition
C4 of Lemma \ref{lem:UPE} holds. Moreover, according to  A1 in Assumption \ref{assume:nonlinearities}, the restriction of
$\alphavec$ on $\Real_{\geq t_0}\times\Omega_\lambda$ is $\lambda$-UPE with $T,\mu$. This,
together with (\ref{eq:lambda_condition:1}), implies that
requirements C1--C4 of Lemma \ref{lem:UPE} are satisfied.

Let $\gamma\in(0,\gamma^\ast]$, where $\gamma^\ast$ is defined
as:
\begin{equation}\label{eq:gamma_ast}
\begin{array}{l}
 \gamma^\ast=\frac{\mu r}{2 D_\varphi
\max\{B_\varphi,M_\varphi,D_\varphi\} T^2}\times \\
\ \ \ \ \ \ \ \ \ \ \ \ \ \ \ \  \left(\sqrt{p}
M_\sigma \max_{i} |\omega_{i}|
\frac{\lambda_{i,\max}-\lambda_{i,\min}}{2}\right)^{-1}.
\end{array}
\end{equation}
This and (\ref{eq:lambda_condition:1}) imply that the
requirement (\ref{eq:lambda_condition}) of the lemma is
satisfied. Hence if $\gamma\in(0,\gamma^\ast]$ then according
to  Lemma \ref{lem:UPE} the matrix $\Phi(t,t_0)$ in
(\ref{eq:error_dynamics:explicit}) satisfies
(\ref{eq:exponential_UPE}). This ensures the existence  of
$\rho,D_\rho\in\Real_{>0}$  such that, along
the solutions of (\ref{eq:error_sys}),
(\ref{eq:nonlinear_par_observer})--(\ref{eq:rational_independence}),
for all $t\geq t_0$ we have: $\norm{\bfe(t)} \leq e^{-\rho(t-t_0)}
D_\rho \norm{\bfe(t_0)}
    + D_\rho \int_{t_0}^t e^{-\rho(t-\tau)}\cdot \norm{\bfv(\tau,\hat \lambdavec(\tau),\lambdavec, y(\tau))}{\rm
    d}\tau$.
The functions $\varphivec$, $\bfg$ in the definition of the
function $\bfv$, (\ref{eq:defn_v}), are Lipschitz with respect
to $\lambdavec$. Furthermore, according to (\ref{eq:xi_bound})
the following holds: $\norm{\xivec(t)}\leq \Delta_\xi$.
Therefore, in accordance with (\ref{eq:defn_v}),
(\ref{eq:xi_bound}), and Assumption \ref{assume:varphi}
\begin{eqnarray}\label{eq:v_Lipschitz}
    & &\norm{\bfv(t,\hat \lambdavec(t),\lambdavec,y(t))} \leq D_v \norm{\hat \lambdavec(\tau) - \lambdavec}_{\infty,[t_0,t]} +
    \Delta_\xi, \\
    & &D_v=D_\varphi \norm{\bfB}    \norm{\thetavec} + D_g, \nonumber
\end{eqnarray}
and hence
\begin{eqnarray}\label{eq:qconv_est}
&&    \begin{array}{l}\norm{\bfe(t)} \leq e^{-\rho(t-t_0)} D_\rho \norm{\bfe(t_0)} + \frac{D_\rho D_v}{\rho}\norm{\hat \lambdavec(\tau) - \lambdavec}_{\infty,[t_0,t]}\end{array} \nonumber\\
&&  \begin{array}{l}+ \frac{D_\rho
    \Delta_{\xi}}{\rho}. \end{array}
\end{eqnarray}
Let $j\in\{1,\dots,p\}$, and consider the solution of
\begin{eqnarray}\label{eq:generator:Poisson:auxiliary}
    \dot q_{2j-1} &=& \omega_j\cdot (q_{2j-1} - q_{2j} - q_{2j-1}(q^2_{2j-1}+q^2_{2j})),  \\
    \dot q_{2j} &=& \omega_j \cdot (q_{2j-1} + q_{2j} - q_{2j}(q^2_{2j-1}+q^2_{2j})),\nonumber
\end{eqnarray}
satisfying the initial condition $q_{2j-1}(t_0)=s_{2j-1}(t_0)$,
$q_{2j}(t_0)=s_{2j}(t_0)$; the parameters $\omega_j$ and values
of $s_{2j-1}(t_0)$, $s_{2j}(t_0)$ are supposed to coincide with
those defined in (\ref{eq:nonlinear_par_observer}),
(\ref{eq:initial_conditions}).  The solution of
(\ref{eq:generator:Poisson:auxiliary}) satisfying initial condition
(\ref{eq:initial_conditions}) is obviously unique, and can be
expressed as a function $\bfq:\Real\rightarrow\Real^{2p}$:
$q_{2j-1}(t)=\cos(\omega_j t + \vartheta_j), \
q_{2j}=\sin(\omega_j t+\vartheta_j), \  \vartheta_j\in\Real,
j=\{1,\dots,p\}$. Parameters $\vartheta_j$ are determined in
accordance with: $\cos(\omega_j t_0+\vartheta_j)=s_{2j-1}(t_0),
\ \sin(\omega_j t_0+\vartheta_j)=s_{2j}(t_0)$.   Given that
$\omega_j$ in (\ref{eq:generator:Poisson:auxiliary}) are
rationally-independent we conclude that the $\omega$-limit set
of $(q_1(t,\bfs_0),q_3(t,\bfs_0),\dots,q_{2p-1}(t,\bfs_0))$ is
the hypercube $[-1,1]^p$ (see e.g. \cite{Katok:1999},
Proposition 1.4.1). Consider the function
$\betavec:\Real^{2p}\rightarrow\Real^p$:
\begin{equation}\label{eq:generator:Poisson:auxiliary:beta}
 \begin{array}{l} \beta_j(\bfq)= \lambda_{j,\min}+\frac{\lambda_{j,\max}-\lambda_{j,\min}}{2}({q}_{2j-1}+1),\end{array}
\end{equation}
and define $\bar{\lambdavec}(t)=\betavec(\bfq(t,\bfs_0))$.
System (\ref{eq:generator:Poisson:auxiliary}),
(\ref{eq:generator:Poisson:auxiliary:beta}) satisfies
conditions P1--P3, and hence we can conclude that trajectory
$\bar{\lambdavec}(\cdot)$ satisfies the recurrence property
(\ref{eq:dense_recurrent}):
\begin{equation}\label{eq:dense_recurrent:2}
\begin{split}
 &\forall \ \lambdavec\in\Omega_\lambda, \ \Delta_{\lambda}\in\Real_{>0}, \ t\geq t_0  \\ & \ \ \ \ \ \ \ \ \ \ \  \ \ \ \  \  \exists \ t'>t: \ \|\lambdavec-\bar{\lambdavec}(t')\|<\Delta_\lambda.
 \end{split}
\end{equation}
Noticing that $\bfs(t,\bfs_0)=\bfq(T(t),\bfs_0)$, where
$T(t)=t_0+\gamma \int \nolimits_{t_0}^t
\sigma(\|y(\tau)-\hat{y}(\tau)\|_\varepsilon) {\rm d}\tau$, we
can conclude that for all $t\geq t_0$ the variable $\hat
\lambdavec(t)$ defined in \eqref{eq:nonlinear_par_observer} can
be expressed in terms of $\bar \lambdavec(T(t))$ as
\beq\label{eq:est_lambda}
 \begin{array}{l}\hat\lambdavec(t) = \bar \lambdavec(t_0 + \gamma \int \nolimits_{t_0}^t \sigma(\|y(\tau)-\hat{y}(\tau)\|_\varepsilon) {\rm d}\tau).\end{array}
\eeq
 Denoting $h(t) = t^\prime - t_0 - \gamma \int
\nolimits_{t_0}^t \sigma(\|y(\tau)-\hat{y}(\tau)\|_\varepsilon
) {\rm d}\tau$, where the value of $t'$ is chosen such that
(\ref{eq:dense_recurrent:2}) holds, we arrive at the following
estimate:
\begin{eqnarray}\label{eq:lambda_approx}
\norm{\lambdavec-\hat \lambdavec(t)} &\leq& \norm{\lambdavec-\bar \lambdavec(t')}+\norm{\bar\lambdavec(t')-\hat \lambdavec(t)}\\
&=&\norm{\lambdavec-\bar \lambdavec(t')}+\norm{\bar\lambdavec(t')-\bar \lambdavec(t'-h(t))}.\nonumber
\end{eqnarray}
The function $\bar{\lambdavec}(\cdot)$ is Lipschitz:
$\|\bar{\lambdavec}(t_1)-\bar{\lambdavec}(t_2)\|\leq \sqrt{p}
\max_{i}\frac{|\omega_i|(\lambda_{i,\max}-\lambda_{i,\min})}{2}|t_1-t_2|,
\ t_1,t_2\in\Real_{\geq t_0}$.
Thus (\ref{eq:est_lambda}), (\ref{eq:lambda_approx}) imply that
\begin{eqnarray}\label{eq:lambda_est}
& &    \norm{\lambdavec-\hat \lambdavec(t)} \leq \norm{\lambdavec-\bar \lambdavec(t')}+\norm{\bar\lambdavec(t')-\hat \lambdavec(t)}\leq \Delta_\lambda \nonumber \\
& &+\begin{array}{l} D_\lambda \norm{h(t)}, \ D_\lambda=\sqrt{p} \max_{i}\frac{|\omega_i|(\lambda_{i,\max}-\lambda_{i,\min})}{2}.\end{array}
\end{eqnarray}
Taking \eqref{eq:qconv_est} and \eqref{eq:lambda_est} into
account we can conclude that the dynamics of the combined
system (\ref{eq:neural_model}), (\ref{eq:linear_par_observer}),
(\ref{eq:nonlinear_par_observer})--(\ref{eq:rational_independence})
obeys the following set of constraints:
\begin{eqnarray}\label{eq:close_loop_system}
   \norm{\bfe(t)} &\leq& e^{-\rho(t-t_0)} D_\rho \norm{\bfe(t_0)} + \begin{array}{l}\frac{D_\rho D_v D_\lambda}{\rho}\norm{h(\tau)}_{\infty,[t_0,t]}\end{array} \nonumber \\
    & & \begin{array}{l} + \frac{D_\rho D_v \Delta_\lambda}{\rho}+\frac{D_{\rho}\Delta_\xi}{\rho},\end{array} \\
   h(t)&=&h(t_0)-\begin{array}{l} \gamma \int_{t_0}^t \sigma(\|\bfC^T\bfe_1(\tau)\|_\varepsilon ) {\rm d}\tau. \end{array} \nonumber
\end{eqnarray}
To proceed further we need an auxiliary result below.
\begin{lem}\label{lem:bound}
Consider a system of which the dynamics for all $t\geq t_0$ satisfy the following inequalities
\begin{eqnarray} \label{eq:interconnection}
        & &\norm{\bfx(t)} \leq \varrho(t-t_0)\norm{\bfx(t_0)}+c \norm{h(\tau)}_{\infty,[t_0,t]} + \Delta,\\
       & & \begin{array}{l} - \int_{t_0}^t \gamma_0(\norm{\bfx(\tau)+\bfd(\tau)}_\varepsilon)d\tau \leq h(t) - h(t_0) \leq 0, \end{array} \nonumber
 \end{eqnarray}
where $\bfx:\Real_{\geq t_0}\rightarrow\Real^n$, $h:\Real\rightarrow\Real$
are trajectories reflecting the evolution of the system's
state, $\bfd:\Real\rightarrow\Real^n$,
$\|\bfd(\tau)\|_{\infty,[t_0,\infty)}\leq \Delta_d$ is a
continuous and bounded function on $[t_0,\infty)$,
$\varrho$ is a strictly monotonically decreasing function
with, $\varrho(0)\geq 1$, $\lim_{s\ra\infty}\varrho(s)=0$;
$c,\Delta\in \Real_{>0}$, and
$\gamma_0:\Real\rightarrow\Real_{\geq 0}$:
\begin{equation}\label{eq:gamma_0}
|\gamma_0(s)|\leq D_\gamma |s|, \ D_\gamma\in\Real_{>0}.
\end{equation}
Then  $\bfx(\cdot)$, $h(\cdot)$ in
(\ref{eq:interconnection}) are globally bounded in forward time, for
$t\geq t_0$, provided that the following conditions hold for
some $d\in(0,1)$, $\kappa\in(1,\infty)$:
\begin{eqnarray}
       & & \begin{array}{l} \varepsilon \geq  \Delta \left(1+\varrho(0)\frac{\kappa}{\kappa - d}\right)+\Delta_d \end{array}, \label{eq:epsilon_choice}\\
      & & \begin{array}{l} D_\gamma \leq \frac{\kappa-1}{\kappa} \left[\varrho^{-1}\left(\frac{d}{\kappa}\right)\right]^{-1}
  \frac{h(t_0) }{ \varrho(0)\|\bfx(t_0)\|+ |h(t_0)|c\left(1+ \frac{\kappa {\varrho}(0)}{(1-d)}\right)}.\end{array} \label{eq:gamma_choice}
\end{eqnarray}
\end{lem}
The proof of Lemma \ref{lem:bound} is provided in the Appendix.
Notice that $h(t)$ in (\ref{eq:close_loop_system}) satisfies
$-\gamma D_\sigma \int_{t_0}^t\|\bfe(\tau)\|_\varepsilon d
\tau\leq h(t)-h(t_0)\leq 0$. Indeed,
$\|\bfC^{T}\bfe_1\|_{\varepsilon}\leq \|\bfe\|_{\varepsilon}$
by virtue of definition of $\|\cdot\|_{\varepsilon}$ and
$\bfC$, and the function $\sigma$ in
(\ref{eq:close_loop_system}) is Lipschitz (see
(\ref{eq:generator_poisson:3})). Thus
(\ref{eq:close_loop_system}) is of the form
(\ref{eq:interconnection}) where
\begin{equation}\label{eq:Delta_v_choice}
\begin{split}
&c={D_\rho D_v D_\lambda}/{\rho}, \ \Delta={D_\rho D_v \Delta_\lambda}/{\rho}+{D_{\rho}\Delta_\xi}/{\rho}, \\
 &{\varrho}(s)=D_{\rho} e^{-\rho s}.
\end{split}
\end{equation}
Notice also that because (\ref{eq:dense_recurrent:2}) holds,
the value of $t'$ in (\ref{eq:lambda_est}) can be chosen
arbitrarily large. This implies that the value of $h(t_0)$ in
(\ref{eq:close_loop_system}) may be chosen arbitrarily large
too. Having this in mind, and invoking  Lemmas \ref{lem:UPE},
\ref{lem:bound} we can conclude that choosing
 $\varepsilon,\gamma$ in
(\ref{eq:nonlinear_par_observer}) as
\begin{eqnarray}\label{eq:epsilon_gamma_choice}
&& \begin{array}{l} \varepsilon\geq r_0(\Delta), \ r_0(\Delta)=\Delta \left(1+D_{\rho}\frac{\kappa}{\kappa - d}\right), \end{array} \nonumber\\
&& \begin{array}{l}  0< \gamma < \bar{\gamma}=\min\{\gamma^\ast, D_{\gamma,\infty}\}\end{array}\nonumber\\
&& \begin{array}{l} D_{\gamma,\infty}= \frac{\kappa-1}{D_\sigma \kappa} \left[\ln\left(D_\rho\frac{ \kappa}{d}\right)\right]^{-1} \frac{\rho}{ c(1+ \kappa D_\rho/(1-d))}, \end{array}
\end{eqnarray}
where $\gamma^\ast$ is defined as in (\ref{eq:gamma_ast}),
ensures that the function $h(\cdot)$ in
(\ref{eq:close_loop_system}) is bounded. Given that  $h(\cdot)$  by construction is
monotone and bounded, the Bolzano-Weierstrass
theorem implies that $h(t)$ converges to a limit, and hence
\begin{equation}\label{eq:convervence_lambda_h}
\begin{split}
&\begin{array}{l}\lim_{t\rightarrow\infty}\int_{t_0}^t \sigma(\|\bfC^{T}\bfe_1(\tau)\|_\varepsilon)d\tau=\bar{h}, \ \bar{h}\in\Real,\end{array}\\
&\begin{array}{l}\lim_{t\rightarrow\infty} \hat{\lambdavec}(t)=\lambdavec^\ast, \ \lambdavec^\ast\in\Omega_{\lambda}.\end{array}
\end{split}
\end{equation}
Noticing that $\sigma(\|\bfC^{T}\bfe(\tau)\|_\varepsilon)$ is
uniformly continuous and using Barbalat's lemma we conclude
that
\begin{equation}\label{eq:y_limit}
\begin{array}{l}\lim_{t\rightarrow\infty}\sup_{\tau\in[t,\infty)}\|\bfC^{T}\bfe_1(\tau)\|\leq
\varepsilon.\end{array}
\end{equation}

{\it Part 3.}  Let us rewrite the equation for $\dot{\bfe}_1$ in (\ref{eq:error_sys}) as: \beq\label{eq:error_sys_3}
\begin{split}
\dot{\bfe}_1&= (\bfA+\bfl \bfC^T)\bfe_1 + \bfv_1(t,\hat\thetavec(t),\hat\lambdavec(t), \thetavec,\lambdavec)\\
&+ \bfv_2(t,\hat\lambdavec(t),\lambdavec)+\bfv_3(t),
\end{split}
\eeq
where $\bfv_3(t)=- \xivec(t)$ and
\begin{eqnarray}\label{eq:defn_v_1}
    & &\bfv_1(t,\hat\thetavec,\hat\lambdavec, \thetavec,\lambdavec) = \bfB(\varphivec^T(t,\hat \lambdavec,y(t))\hat{\thetavec}-\varphivec^T(t,\lambdavec,y(t))\thetavec) \nonumber \\
& &  \bfv_2(t,\hat\lambdavec,\lambdavec)=\bfg(t,\hat{\lambdavec},y(t),u(t))-\bfg(t,\lambdavec,y(t),u(t)). \label{eq:defn_v_2}
    \end{eqnarray}
Next steps make use of the following lemma. {
\begin{lem}\label{lem:observer_inferrence}
Consider
\begin{equation}\label{eq:system_io}
\begin{split}
& \begin{array}{ll}
\dot{\bfx}&=\bfA \bfx + \bfu(t)+\bfd(t),\\
y&=\bfC^{T}\bfx, \ \bfx(t_0)=\bfx_0, \ \bfx_0\in\Real^n,
\end{array}
\end{split}
\end{equation}
where $\bfA$ and $\bfC$ are defined as in
(\ref{eq:observer_canonical_form}), and
$\bfu,\bfd:\Real\rightarrow\Real^n$, $\bfu\in\mathcal{C}^1$,
$\bfd\in\mathcal{C}$. Let $\bfu,\dot{\bfu},\bfd$ be bounded:
$\max\{\|\bfu(t)\|,\|\dot{\bfu}(t)\|\}\leq B, \
\|{\bfd}(t)\|\leq \Delta_{\xi}$ for all $t\geq t_0$.

Then, if the solution of (\ref{eq:system_io}) is globally bounded for all $t\geq t_0$, there exist $\kappa_1,\kappa_2\in\mathcal{K}$:
\[
\begin{split}
&\|y(\tau)\|_{\infty,[t_0,\infty)}\leq \varepsilon \Rightarrow  \ \exists \ t'(\varepsilon,\bfx_0)\geq t_0: \\
&\ \ \ \ \ \ \ \left\|z_1(\tau) + u_1(\tau) \right\|_{\infty,[t',\infty)}\leq \kappa_1(\varepsilon) + \kappa_2(\Delta_{\xi}),
\end{split}
\]
where $z_1=(1 \ 0 \ \dots \ 0)\bfz$,
\begin{equation}\label{eq:system_io_filter}
\begin{split}
\dot\bfz&= \Lambda \bfz +  \bfG \bfu(t),\ \Lambda=\left(\begin{array}{c|c}
\multirow{2}{*}{\ensuremath{-\bfb}} &   I_{n-2}\\ & 0 \end{array}\right),\\
\bfG&=\left(\begin{array}{cc} - {\bfb} & I_{n-1} \end{array}\right), \ \bfz(t_0)=0,
\end{split}
\end{equation}
and  $\bfb=\col{b_1,\dots,b_{n-1}}$:  real parts of the roots
of $s^{n-1}  + b_1 s^{n-2}+\cdots + b_{n-1}$ are negative.

Moreover, if $\bfd(t)\equiv0$, then
\begin{equation}\label{eq:identifiability}
\begin{split}
&y(t)=0 \ \forall \ t\geq t_0 \Rightarrow \ \exists \ \bfp\in\Real^{n-1}: \ \forall \ t\geq t_0\\
&(1 \ 0 \ \dots \ 0)e^{\Lambda(t-t_0)}\bfp+z_1(t)+u_1(t)=0 .
\end{split}
\end{equation}
\end{lem}
} The proof of Lemma \ref{lem:observer_inferrence} is provided
in the Appendix. According to (\ref{eq:xi_bound}) and
Assumption \ref{assume:varphi},  $\bfv_1(\cdot,\hat\thetavec(\cdot),\hat\lambdavec(\cdot),\thetavec,\lambdavec)$,
$\bfv_2(\cdot,\hat\lambdavec(\cdot),\lambdavec)$, $\bfv_3(\cdot)$ and $\dot{\bfv}_1$,
$\dot{\bfv}_2$ are bounded. Thus assumptions of Lemma
\ref{lem:observer_inferrence} are met for equations
(\ref{eq:error_sys_3}), (\ref{eq:defn_v_1}), and hence
(\ref{eq:y_limit}) implies that there is a
$t_1(\varepsilon)\geq t_0$ and
$\kappa_1,\kappa_2\in\mathcal{K}$ such that $\forall \ t \geq
t_1(\varepsilon)$ we have:
\[
\begin{array}{l}
\big\|\varphivec^T(t,\hat \lambdavec(t),y(t))\hat{\thetavec}(t)-\varphivec^T(t,\lambdavec,y(t))\thetavec +v_{2,1}(t,\hat\lambdavec(t),\lambdavec)\\
+\tilde{\bfC}^{T}\int_{t_0}^{t} e^{\Lambda(t-\tau)}\bfG \bfv_2(\tau,\hat\lambdavec(\tau),\lambdavec) d\tau\big\|  \leq\kappa_1(\varepsilon)+\kappa_2(\Delta_\xi),
\end{array}
\]
where $\bfG\in\Real^{(n-1)\times n}$,
$\tilde{\bfC}\in\Real^{n-1}$ are defined as in
(\ref{eq:goal_lambda})--(\ref{eq:unidentifable_set_2}), and $v_{2,1}(\cdot)$ is the first component of $\bfv_2(\cdot)$. Given
that $\int_{t_0}^{t} e^{\Lambda(t-\tau)}\bfG
\bfv_2(\tau,\hat\lambdavec(\tau),\lambdavec) d\tau =
\int_{t_0}^{t} e^{\Lambda(t-\tau)}\bfG
\bfv_2(\tau,\lambdavec^\ast,\lambdavec) d\tau$ $+$
$\int_{t_0}^{t} e^{\Lambda(t-\tau)}\bfG
(\bfv_2(\tau,\hat\lambdavec(\tau),\lambdavec)$ $-
\bfv_2(\tau,\lambdavec^\ast,\lambdavec)) d\tau$, noticing that
$\Lambda$ is Hurwitz and  that, according to
(\ref{eq:convervence_lambda_h})
$\bfv_2(t,\hat\lambdavec(t),\lambdavec)-
\bfv_2(t,\lambdavec^\ast,\lambdavec)\rightarrow 0$ as
$t\rightarrow\infty$,  we can conclude that there is a
$t_2(\varepsilon)\geq t_1(\varepsilon)$ such that
$\etavec(t,\hat{\thetavec}(t),\lambdavec^\ast,\thetavec,\lambdavec)$
defined as in (\ref{eq:unidentifable_set_1}) satisfies:
\begin{eqnarray}\label{eq:nonlinear_PE_regressor_small}
&&\begin{array}{l}\|\etavec(t,\hat{\thetavec}(t),\lambdavec^\ast,\thetavec,\lambdavec)\|=\big\|\varphivec^T(t, \lambdavec^\ast,y(t))\hat{\thetavec}(t)-\end{array}\nonumber\\
&&\begin{array}{l}\varphivec^T(t,\lambdavec,y(t))\thetavec +\tilde{\bfC}^{T}\int_{t_0}^{t} e^{\Lambda(t-\tau)}\bfG \bfv_2(\tau,\lambdavec^\ast,\lambdavec) d\tau\end{array}\nonumber\\
&&\begin{array}{l}+v_{2,1}(t,\lambdavec^\ast,\lambdavec)\big\| \leq\kappa_1(\varepsilon)+\kappa_2(\Delta_\xi)+\varepsilon \ \ \forall  \ t\geq t_2(\varepsilon).\end{array}
\end{eqnarray}
Recall that the restriction of
$\alphavec_2(t,(\lambdavec',\thetavec'),(\lambdavec,\thetavec))=\etavec(t,{\thetavec}',\lambdavec',\thetavec,\lambdavec)$ on $\Real_{\geq t_0}\times\Real^{p+m}\times\Real^{p+m}$ is wNPE with $L,\beta,\mathcal{E}$.
Let $t_3(\varepsilon)$ be such that
$\|\bfC^{T}\bfe_1(t)\|<2\varepsilon$ for all $t\geq
t_3(\varepsilon)$ (existence of such $t_3(\varepsilon)$ follows
from (\ref{eq:y_limit})). Consider the sequence
$\{\tau_i\}_{i=0}^{\infty}$,
$\tau_i=\max\{t_3(\varepsilon),t_2(\varepsilon)\}+i L$. Since $\varphivec(\cdot,\hat{\lambdavec}(\cdot),y(\cdot))$ is bounded, there
is an $M_{\theta}\in\Real_{>0}$:
\begin{equation}\label{eq:nonlinear_PE_regressor_small:1}
\|\hat{\thetavec}(\tau)-\hat{\thetavec}(\tau_i)\|_{\infty,[t_i,t_{i+1}]}\leq \varepsilon 2 \gamma_\theta B_\varphi L =\varepsilon M_{\theta}
\end{equation}
for all $t\geq \tau_0$.  Hence $\forall \ t: \ t\in[\tau_i,\tau_{i+1}], \ i\in\Natural$, we have:
$\|\etavec(t,\hat{\thetavec}(\tau_i),\lambdavec^\ast,\thetavec,\lambdavec)\|\leq \kappa_1(\varepsilon)+\kappa_2(\Delta_\xi)+\varepsilon(M_{\theta} B_\varphi+1)$.
This, however, implies that there is an $N\in\Natural$ such that $\dist\left(\left(\begin{array}{c}\lambdavec^\ast\\
                                \hat\thetavec(\tau_i)
    \end{array}\right),\mathcal{E}(\lambdavec,\thetavec)\right)\leq \beta^{-1}(\kappa_1(\varepsilon)+\kappa_2(\Delta_\xi)+\varepsilon(M_{\theta} B_\varphi+1))$
for all $i\geq N$. Therefore, taking
(\ref{eq:lambda_condition:1}),
(\ref{eq:nonlinear_PE_regressor_small:1}) into account, we can
conclude that there is a $t'(\varepsilon)$:
\[
\begin{split}
& \dist\left(\left(\begin{array}{c}\hat\lambdavec(t)\\
                                \hat\thetavec(t)
    \end{array}\right),\mathcal{E}(\lambdavec,\thetavec)\right)\leq 2\varepsilon M_{\theta}+\\
& \ \ \ \ \beta^{-1}(\kappa_1(\varepsilon)+\kappa_2(\Delta_\xi)+\varepsilon(M_{\theta} B_\varphi+1))  \ \forall \ t\geq t'(\varepsilon).
\end{split}
\]
Notice that $r_0$ in (\ref{eq:epsilon_gamma_choice}) is a class
$\mathcal{K}_{\infty}$ function of $\Delta$.
Parameter $\Delta$, as defined in (\ref{eq:Delta_v_choice}),
is the sum: $\Delta=\frac{D_\rho D_v
\Delta_\lambda}{\rho}+\frac{D_{\rho}\Delta_\xi}{\rho}$. Given
that the value of $\Delta_\lambda$ can be chosen arbitrarily,
we  pick $\Delta_\lambda=\Delta_\xi$. This renders  $r_0$ in (\ref{eq:epsilon_gamma_choice}) a class
$\mathcal{K}_\infty$ (and hence class $\mathcal{K}$) function
of  $\Delta_\xi$. Denote this function as $r_1$, then $\varepsilon > r_1(\Delta_\xi)$ implies that
\begin{eqnarray}\label{eq:r_1_2_estimate}
&&\beta^{-1}(\kappa_1(\varepsilon)+\kappa_2(\Delta_\xi)+\varepsilon(M_{\theta} B_\varphi+1))+2\varepsilon
M_{\theta} <\\
&&
\beta^{-1}(\kappa_1(\varepsilon)+\kappa_2(r_1^{-1}(\varepsilon))+\varepsilon(M_{\theta} B_\varphi+1))\nonumber\\
& & +2\varepsilon
M_{\theta}=r_2(\varepsilon).\nonumber
\end{eqnarray}
Thus (\ref{eq:thm2}) holds.

Finally, if $\mathcal{E}(\lambdavec,\thetavec)$ coincides with
$\mathcal{E}_0(\lambdavec,\thetavec)$, then Assumption
\ref{assume:varphi} and (\ref{eq:thm2}) imply that
$\|\bfv_1(t,\hat\thetavec(t),\hat\lambdavec(t), \thetavec,\lambdavec)
+ \bfv_2(t,\hat\lambdavec(t),\lambdavec)\|<M_1 r_2(\varepsilon)$
for some $M_1\in\Real_{>0}$, $t\geq t'(\varepsilon)$. Since $\bfA+\bfl\bfC^T$ in
(\ref{eq:error_sys_3}) is Hurwitz,  (\ref{eq:thm1})
follows. $\square$


\section{Discussion and generalization}\label{sec:discussion}

\subsection{Removing passivity requirement (Assumption \ref{assume:matrices})}

Theorem \ref{theorem:neural_identification} requires that
$\bfA$, $\bfB$,  $\bfC$ in (\ref{eq:neural_model}) satisfy
Assumption \ref{assume:matrices}. Here  we invoke the idea of filtered transformations \cite{Marino92}, \cite{MarinoTomei93} to show how observer
(\ref{eq:linear_par_observer}),
(\ref{eq:nonlinear_par_observer}) can be modified so that this
condition could be replaced with the requirement that the pair $\bfA$, $\bfC$ is observable. Consider a generalization of (\ref{eq:neural_model})
\begin{eqnarray}\label{eq:neural_model:3}
\dot{\bfx}&=&\bfA \bfx + \Psivec(t,\lambdavec,y)\thetavec+\bfg(t,\lambdavec,y,u(t))+\xivec(t),\\
y&=&\bfC^{T}\bfx, \ \bfA=\left(\begin{array}{cc}0 & \bfI_{n-1}\\ 0 &
0\end{array}\right), \ \bfC=\mathrm{col}(1,0,\dots,0), \nonumber
\end{eqnarray}
where
$\Psivec:\Real\times\Real^p\times\Real\rightarrow\Real^{n\times
m}$, $\Psivec\in\mathcal{C}^1$, is  Lipschitz in $\lambdavec$,
 and ${\Psivec}(\cdot,\lambdavec,y(\cdot))$,
$\dot{\Psivec}(\cdot,\lambdavec,y(\cdot))$ are bounded on $\Real_{\geq t_0}$.
 The function $\xivec$ and parameters are
defined as in (\ref{eq:neural_model}), and the function $\bfg$
satisfies Assumption \ref{assume:varphi}.

Let $\bfB=\col{1,b_1,\dots,b_{n-1}}$ be a vector such that the
polynomial $s^{n-1}+b_1 s^{n-2}+\dots+b_{n-1}$ is Hurwitz. As
an observer candidate for (\ref{eq:neural_model:3}) we propose
a system in which $\mathcal{S}_w$ is defined as in
(\ref{eq:nonlinear_par_observer}), and $\mathcal{S}_a$ is given
as follows:
\begin{eqnarray}\label{eq:linear_par_observer:modified}
&&\dot{\bfM}=(\bfA-\bfB\bfC^{T}\bfA)\bfM+(\bfI_n-\bfB\bfC^{T})\Psivec(t,\hat{\lambdavec}(t),y(t)),\nonumber\\
&&\dot{\hat{\zetavec}}=\bfA\hat{\zetavec}+\bfl(\bfC^{T}\hat{\zetavec}-y(t))+\bfB{\varphivec}^{T}(t,\hat{\lambdavec}(t),y(t),[\hat\lambdavec,y])\hat{\thetavec}\nonumber\\
&& \ \ \ \ \ \ \ \ \ \ \ \ \ \ + \bfg(t,\hat{\lambdavec}(t),y(t),u(t)),  \\
&&\dot{\hat{\thetavec}}=-\gamma_\theta (\bfC^{T}\hat{\zetavec}-y(t)){\varphivec}(t,\hat{\lambdavec}(t),y(t),[\hat\lambdavec,y]), \ \gamma_\theta\in\Real_{>0},\nonumber\\
&& \hat{\bfx}=\hat{\zetavec}+\bfM\hat{\thetavec}, \ \ \bfM\in\Real^{n\times m}, \ \ \bfM(t_0)=0\nonumber,
\end{eqnarray}
where
\begin{equation}\label{eq:regressor_modified}
\begin{split}
{\varphivec}^{T}(t,\hat{\lambdavec}(t),y(t),[\hat\lambdavec,y]) =&\bfC^{T}\bfA \bfM(t,[\hat{\lambdavec},y])\\ &+\bfC^{T}\Psivec(t,\hat{\lambdavec}(t),y(t)).
\end{split}
\end{equation}
The first row  of $\bfM$ is zero for all $t\geq t_0$, and $\hat{y}=\bfC^{T}\hat{\bfx}=\bfC^{T}\hat{\zetavec}$.
Since  $\Psivec(\cdot,\hat{\lambdavec}(\cdot),y(\cdot))$ is bounded on $\Real_{\geq t_0}$ and
Lipschitz in $\hat{\lambdavec}$,
$\bfM(\cdot,[\hat{\lambdavec},y]), \dot{\bfM}(\cdot,[\hat{\lambdavec},y])$
are globally bounded on $\Real_{\geq t_0}$, and ${\bfM}(t,[\hat\lambdavec,y])$ is Lipschitz
in $\hat{\lambdavec}$ for $\hat\lambdavec=\mathrm{const}$. Let
$\zetavec=\bfx-\bfM\thetavec$, then  using (\ref{eq:neural_model:3})--(\ref{eq:regressor_modified}) we can write
\[
\begin{split}
{\dot\zetavec}=&\bfA \zetavec
+ \bfB\varphivec(t,\hat\lambdavec(t),y(t),[\hat{\lambdavec},y])\thetavec +
(\Psivec(t,\lambdavec,y(t),u(t))\\
&-\Psivec(t,\hat\lambdavec(t),y(t),u(t)))\thetavec
+ \bfg(t,\lambdavec,y(t),u(t))+\xivec(t).
\end{split}
\]
Dynamics of (\ref{eq:neural_model:3}),
(\ref{eq:linear_par_observer:modified}) in the coordinates
$\bfe_1=\hat{\zetavec}-\bfx+\bfM\thetavec$,
$\bfe_2=\hat{\thetavec}-\thetavec$ is
\begin{eqnarray}\label{eq:error_sys:modified}
&& \bpm \dot \bfe_1 \\ \dot \bfe_2 \epm
        =
        \bpm
            \bfA+\bfl \bfC^T   &   \bfB \alphavec^{T}(t) \\
            -\gamma_{\theta} \alphavec(t)\bfC^T & 0
        \epm  \bpm
            \bfe_1 \\ \bfe_2
        \epm \nonumber\\
&& \ \ \ \ \ \ \ \ \ \ \ \ \ \ \
         +
        \bpm
            \tilde{\bfv}(t,\hat \lambdavec(t),\lambdavec, y(t),u(t)) \\ 0
        \epm
\end{eqnarray}
where $\alphavec(t)=\varphivec(t,\hat \lambdavec(t),y(t),[\hat \lambdavec,y])$, $\tilde{\bfv}(t,\hat \lambdavec,\lambdavec,
y,u)=(\Psivec(t,\hat
\lambdavec,y)-\Psivec(t,{\lambdavec},y))\thetavec +\bfg(t,\hat
\lambdavec,y,u)-\bfg(t,\lambdavec,y,u)-\xivec(t)$. Since the
pair $\bfA$, $\bfC$ is observable there always is an $\bfl$
so that (\ref{eq:condition_MKY}) holds.  The structure of
(\ref{eq:error_sys:modified}) is now identical to that of
(\ref{eq:error_sys}), and Assumptions
\ref{assume:matrices}, \ref{assume:varphi} hold for the
functions $\varphivec$, $\bfg$ in
(\ref{eq:linear_par_observer:modified}). Finally, consider the
function $\etavec_1$:
\[
\begin{split}
&\etavec_1(t,\lambdavec,\thetavec,\lambdavec',\thetavec')=\varphivec^{T}(t,\lambdavec',y(t),[\lambdavec',y])(\thetavec'-\thetavec)+\\
&\bfC^{T}(\Psivec(t,\lambdavec',y(t))-\Psivec(t,\lambdavec,y(t)))\thetavec + g_1(t,\lambdavec',y(t),u(t))-\\
& g_1(t,\lambdavec,y(t),u(t))+ q(t,\lambdavec',\lambdavec,\thetavec),
\end{split}
\]
where $q(t,\lambdavec',\lambdavec,\thetavec)=\tilde{\bfC}\bfz(t,\lambdavec,\lambdavec',\thetavec)$,
$\dot{\bfz}=\Lambda \bfz + \bfG(
(\Psi(t,\lambdavec',y(t))-\Psi(t,\lambdavec,y(t)))\thetavec+\bfg(t,\lambdavec',y(t),u(t))-\bfg(t,\lambdavec,y(t),u(t)))$,
$\bfz(t_0)=0$,  and $\tilde{\bfC}$, $\Lambda$, $\bfG$ are
defined as in (\ref{eq:unidentifable_set_2}). The following is
now immediate:
\begin{thm}\label{cor:general} Consider \eqref{eq:neural_model:3},  \eqref{eq:linear_par_observer:modified}, \eqref{eq:nonlinear_par_observer}--\eqref{eq:rational_independence}. Suppose that  condition A1 of Assumption \ref{assume:nonlinearities} holds for the function $\alphavec_3:\Real_{\geq t_0}\times\Omega_\lambda\rightarrow\Real^m$,  $\alphavec_3(t,\lambdavec)=\varphivec(t,\lambdavec,y(t),[\lambdavec,y])$, where $\varphivec$ is defined as in (\ref{eq:linear_par_observer:modified}). Furthermore, let the restriction of
$\alphavec_4:\Real\times\Real^{p+m}\times\Real^{p+m}\rightarrow\Real$, $\alphavec_4(t,(\lambdavec,\thetavec),(\lambdavec',\thetavec'))=\etavec_1(t,\lambdavec,\thetavec,\lambdavec',\thetavec')$ on $\Real_{\geq t_0}\times\Real^{p+m}\times\Real^{p+m}$
be weakly nonlinearly persistently exciting in
$(\lambdavec,\thetavec)$ wrt to the map $\mathcal{E}_1$:
\[
\begin{split}
&\mathcal{E}_1(\lambdavec,\thetavec)=\{(\lambdavec',\thetavec'), \ \lambdavec'\in\Real^p, \ \thetavec'\in\Real^m| \bfB (\thetavec'-\thetavec)^{T}\cdot \\
& \varphivec(t,\lambdavec',y(t),[\lambdavec',y])+ (\Psivec(t,\lambdavec',y(t))-\Psivec(t,\lambdavec,y(t)))\thetavec\\
&+\bfg(t,\lambdavec',y(t),u(t))-\bfg(t,\lambdavec,y(t),u(t))=0,\ \forall \ t\geq t_0\}.
\end{split}
\]
Then there exist a constant $\bar \gamma \in \Real_{>0}$ and
functions $r_1,r_2,r_3\in\mathcal{K}$ such that if $\gamma,
\varepsilon$ are the corresponding parameters of
(\ref{eq:nonlinear_par_observer}), and $\gamma \in (0, \bar
\gamma)$, $\varepsilon > r_1(\Delta_{\xi})$,  then
(\ref{eq:thm1}), (\ref{eq:thm2}) hold (with $\mathcal{E}$ replaced by $\mathcal{E}_1$) for  the interconnection
(\ref{eq:neural_model:3}),
(\ref{eq:linear_par_observer:modified}), and
(\ref{eq:nonlinear_par_observer}).
\end{thm}
The proof is largely identical to that of Theorem
\ref{theorem:neural_identification} (a sketch is presented in
the Appendix).  According to Remarks \ref{rem:relax_pe}, \ref{rem:relax}
one can replace the requirement that the restriction of $\alphavec_4$ on $\Real_{\geq t_0}\times\Real^{p+m}\times\Real^{p+m}$ is {wNPE} with $L,\beta,\mathcal{E}_1$ with that of the $\lambda'$-uniform persistency of excitation of
the the restriction of $\alphavec_5$:
\begin{equation}\label{eq:nonlinear_pe_new}
\alphavec_5(t,\lambdavec')=(\varphivec^{T}(t,\lambdavec',y(t),[\lambdavec',y]),\bfR_1(t,\lambdavec,\lambdavec',\thetavec)),
\end{equation}
where $\bfR_1(t,\lambdavec,\lambdavec',\thetavec)= \int_{0}^1
\frac{\pd }{\pd \bfs}
r_1(t,\bfs(\xi,\lambdavec,\lambdavec'),\thetavec) d\xi$,
$\bfs(\xi,\lambdavec,\lambdavec')=\lambdavec' \xi +
(1-\xi)\lambdavec$, and $r_1(t,\lambdavec,\thetavec)$ $=$
$\bfC^{T}$ $\Psivec(t,\lambdavec,y(t))\thetavec$ $+$
$g_1(t,\lambdavec,y(t),u(t))$ $+\tilde{\bfC}^{T}$ $\int_{t_0}^t
e^{\Lambda(t-\tau)}\bfG$
$(\Psivec(\tau,\lambdavec,y(\tau))\thetavec+\bfg(\tau,\lambdavec,y(\tau),u(\tau)))d\tau$, on $\Real_{\geq t_0}\times\Omega_\lambda$.

Consider now systems (\ref{eq:neural_model:ext}). Since $\bfA,\bfC$ is observable, there is a
coordinate transformation $\bfx\mapsto T(\bfA) \bfx$ bringing  system (\ref{eq:neural_model:ext}) into the form (\ref{eq:neural_model:3}), albeit with the functions $\Psivec$, $\bfg$ and vector $\thetavec$ defined differently. An example illustrating the viability of this approach is  provided in Section \ref{sec:example}.
Notice also that observability of $\bfA$,
$\bfC$ implies that the system $\dot{\bfx}=\bfA \bfx +
\tilde\Psivec(t,\lambdavec,\bfx)\thetavec+\tilde\bfg(t,\lambdavec,\bfx,u(t))+\xivec(t)$,
$y=\bfC^{T}\bfx$, in which the functions $\tilde\Psivec$,
$\tilde\bfg$ are bounded and Lipschitz in $\bfx$ can be brought
into the form (\ref{eq:neural_model:3}) by using an auxiliary
high-gain observer (cf. \cite{SysConLett:Grip:2011}).

\subsection{Presence of measurement noise}

{ Suppose now that observations of system
(\ref{eq:neural_model}) output, $y$, are corrupted by noise.
That is, instead of $y=\bfC^T\bfx$ we can access only the
variable $y_d=\bfC^{T}\bfx + d$,
$y_d\in\mathcal{D}_y$, where $d:\Real\rightarrow\Real$,
$d\in\mathcal{C}^1$, $\|d(\tau)\|_{\infty,[t_0,\infty)}\leq
\Delta_d$, $\Delta_d\in\Real_{\geq 0}$, and $|\dot{d}(t)|$ is
bounded. In this case the variable $y$ in the observer
definition (\ref{eq:linear_par_observer}),
(\ref{eq:nonlinear_par_observer}) is replaced by $y_d$, and the
dynamics of $\bfe_1=\hat{\bfx}-\bfx$,
$\bfe_2=\hat{\thetavec}-\thetavec$ becomes:
\[
    \splt{
        \bpm \dot \bfe_1 \\ \dot \bfe_2 \epm
        =&
        \bpm
            \bfA+\bfl \bfC^T   &   \bfB \alphavec^T(t,\hat \lambdavec(t)) \\
            -\gamma_{\theta} \alphavec(t,\hat \lambdavec(t))\bfC^T & 0
        \epm
        \bpm
            \bfe_1 \\ \bfe_2
        \epm         \\
        & +
        \bpm
            \bfv(t,\hat \lambdavec(t),\lambdavec, y_d(t),u(t)) \\ 0
        \epm+\bpm
            \xivec_1(t) \\ \xivec_2(t)
        \epm
    }
\] where $\alphavec^T(t,\hat
\lambdavec)=\varphivec^{T}(t,\hat{\lambdavec},y_d(t))$,
$\bfv(t,\hat\lambdavec, \lambdavec,y_d,u)  =
\bfB(\varphivec^T(t,\hat
\lambdavec,y_d)-\varphivec^T(t,\lambdavec,y_d))\thetavec
     +  \bfg(t,\hat{\lambdavec},y_d,u)-\bfg(t,\lambdavec,y_d,u) -
     \xivec(t)$, and $\xivec_1(t)=\bfB(\varphivec^T(t,
\lambdavec,y_d(t))-\varphivec^T(t,
\lambdavec,y(t)))\thetavec+(\bfg(t,\lambdavec,y_d(t),u(t))-\bfg(t,\lambdavec,y(t),u(t)))-\bfl d(t)$,
$\xivec_2(t)=-\gamma_\theta d(t) \alphavec(t,\hat \lambdavec(t))$.
It can now be seen that if  $\varphivec$, $\bfg$ are
Lipschitz in $y$ then there is an $M_d>0$ such that
$\max\{\|\xivec_1(\tau)\|_{\infty,[t_0,\infty)},\|\xivec_2(\tau)\|_{\infty,[t_0,\infty)}\}\leq
M_d\Delta_d$. Thus invoking Lemmas \ref{lem:UPE}, \ref{lem:bound}
and following the argument provided in proof of Theorem
\ref{theorem:neural_identification} one can establish existence
of $\gamma>0$ and $\varepsilon>0$ such that
(\ref{eq:convervence_lambda_h}), (\ref{eq:y_limit}) hold.
Convergence of the estimates will also follow subject to corresponding persistency of excitation
requirements (cf. part 3 of the
proof).} An illustration of the influence of measurement noise on performance of the observer is provided in Section \ref{sec:example}, Fig. \ref{fig:example:1}.

\section{Examples}\label{sec:example}

Consider  the following system:
\begin{eqnarray}\label{eq:example_academic:1}
\dot{\bfx}&=&\bfA \bfx +\bfB\theta + \bfB(\sin(\lambda \cos(t))+e^{\lambda \sin(t)})+\xivec(t),  \\
 \bfA&=&\left(\begin{array}{cc}
- 2 & 1\\
-1 & 0
\end{array}\right), \ \bfB=\left(\begin{array}{c}
1\\
1
\end{array}\right), \ \begin{array}{l} y=\bfC^{T}\bfx,  \ \bfx(t_0)=\bfx_0,\\ \bfC=\mathrm{col}(1,0,\dots,0),\end{array} \nonumber
\end{eqnarray}
where $\theta\in [0,1]=\Omega_\theta$,
$\lambda\in[0.1,1]=\Omega_\lambda$ are unknown parameters, and
$\bfx_0$ is only partially known. The function
$\xivec:\Real\rightarrow\Real^2$,
$\xivec(t)=0.001\col{\sin(t),\cos(t)}$, stands for the
unmodeled dynamics and is supposed to be unavailable for direct
observation.

Let the task be to infer the values of $\bfx$, $\theta$,
$\lambda$ from the measurements of $y$ over time. System
(\ref{eq:example_academic:1}) belongs to the class of equations
described by (\ref{eq:neural_model}) with
$\varphivec(t,\lambda,y)=1$ $\forall \ t,\lambda,y$, and
$\bfg(t,\lambda,y,u)=\bfB(\sin(\lambda \cos(t))+e^{\lambda
\sin(t)})$. Moreover, it satisfies Assumption
\ref{assume:matrices} with
\[
\bfP=\left(\begin{array}{cc}2 & -1\\ -1 & 1 \end{array}\right), \ \bfQ= \left(\begin{array}{cc}6 & -3\\ -3 & 2 \end{array}\right), \ \bfl=\left(\begin{array}{c}0\\ 0\end{array}\right),
\]
and Assumption \ref{assume:varphi} with $D_\varphi=0$,
$D_g=B_g=M_g=\sqrt{2}(1+e)$, $B_\varphi=1$, and $M_\varphi=0$.
According to (\ref{eq:linear_par_observer})--
(\ref{eq:nonlinear_par_observer}) the  observer candidate is:
\begin{eqnarray}
&&\begin{split}
\dot{\hat\bfx}&=\bfA\hat\bfx+\bfB\hat{\theta}+\bfB(\sin(\hat\lambda \cos(t))+e^{\hat\lambda \sin(t)})\\
\dot{\hat{\theta}}&=-\gamma_\theta (\bfC^{T}\hat{\bfx}-y(t))
\end{split}\label{eq:example_academic:2}\\
&&\dot{s}_{1}=\gamma \tanh(\|\bfC^{T}\hat{\bfx}-y(t)\|_\varepsilon)  (s_{1} - s_{2} - s_{1}(s_{1}^2+s_{2}^2)) \nonumber \\
&& \dot{s}_{2}=\gamma \tanh(\|\bfC^{T}\hat{\bfx}-y(t)\|_\varepsilon)  (s_{1}+ s_{2} - s_{2}(s_{1}^2+s_{2}^2))\nonumber \\
&&\hat{\lambda}= 0.1+0.45(s_{1}+1), \ s_{1}^2(t_0)+s_{2}^2(t_0)=1.\label{eq:example_academic:3}
\end{eqnarray}
Parameters $\gamma$, $\gamma_\theta$, and $\varepsilon$ will be
specified in a later stage.

Note that
$\varphivec(t,\lambda,y)=1$. Hence  condition A1 in Assumption \ref{assume:nonlinearities}
holds. Notice also that the function $\etavec$,
as defined in (\ref{eq:unidentifable_set_1}), in this case
becomes:
$\etavec(t,\lambda,\theta,\lambda',\theta')=\theta-\theta'+\sin(\lambda
\cos(t))+e^{\lambda \sin(t)}-\sin(\lambda' \cos(t))-e^{\lambda'
\sin(t)}$. One can now check that condition  A2 is satisfied as well. According to Theorem
\ref{theorem:neural_identification}, system
(\ref{eq:example_academic:2}), (\ref{eq:example_academic:3}) is
an adaptive observer for (\ref{eq:example_academic:1}) subject
to the choice of $\gamma$, $\gamma_\theta$, and $\varepsilon$.
A procedure for setting specific values of these parameters can
be derived from the argument provided in the proof of the
theorem. Let us show how this procedure works in this
example.

According to the theorem, parameter $\gamma_\theta$ is an
arbitrary positive number; here, for simplicity, we set
$\gamma_\theta=1$. Parameters $\gamma$,
$\varepsilon$ are to satisfy
(\ref{eq:Delta_v_choice}), (\ref{eq:epsilon_gamma_choice}).
The choice of  $\gamma$ is subjected to
two constraints. The first constraint is
$\gamma\in(0,\gamma^\ast]$, where $\gamma^\ast$ is specified in
(\ref{eq:gamma_ast}). It ensures that the restriction of
$\varphivec(\cdot,\hat{\lambda}(\cdot),y(\cdot))$  on $\Real_{\geq t_0}$ is persistently exciting.
In our case
$\varphivec(t,\hat{\lambda}(t),y(t))$ is independent on
$\hat{\lambda}(t)$, and this property holds for any $\gamma>0$. The
second constraint is: $\gamma<\frac{\kappa-1}{D_\sigma \kappa}
\left[\ln\left(D_\rho\frac{ \kappa}{d}\right)\right]^{-1}
\frac{\rho}{ c(1+ \kappa D_\rho/(1-d))}, \ \kappa>1,  \
d\in(0,1), \ D_\sigma=1, \ c={D_\rho D_v D_\lambda}/{\rho}$, where
$D_\lambda=0.45$, $D_v=\sqrt{2}(1+e)$ (see
(\ref{eq:lambda_est}), (\ref{eq:v_Lipschitz})), and $\rho$ and
$D_\rho$ are such that the fundamental matrix of solutions of
(\ref{eq:sys_lemma_exponential}), $\Phi(t,t_0)$, satisfies:
$\|\Phi(t,t_0)\bfp\|\leq D_\rho e^{-\rho(t-t_0)}\|\bfp\|$. In
this example
\[
\Phi(t,t_0)=e^{\bfA_1(t-t_0)}, \ \bfA_1=\left(\begin{array}{ccc}-2 & 1 & 1\\ -1 & 0& 1\\ -1 & 0& 0\end{array}\right),
\]
and $\rho=0.5$, $D_\rho=4.242$. Picking $d=0.2$, $\kappa=2$
results in $\bar\gamma=0.00286$. Since the values of $d$,
$\kappa$ are now defined, we can set the value of
$\varepsilon$. Taking (\ref{eq:Delta_v_choice}),
(\ref{eq:epsilon_gamma_choice}) into account and noticing that
$\Delta_\lambda=0$  (this is because trajectories of
(\ref{eq:example_academic:3}) in which the term
$\tanh(\|\bfC^{T}\hat{\bfx}-y(t)\|_\varepsilon)$ is replaced with
$1$ will pass through every point of
$\Omega_\lambda=[0.1,1]$), we arrive at $\varepsilon\geq
\frac{D_\rho
\Delta_\xi}{\rho}\left(1+D_\rho\frac{\kappa}{\kappa-d}\right)$,
where $\Delta_\xi: \ \|\xivec(t)\|\leq \Delta_\xi$. Given that
$\Delta_\xi=0.001$ we obtain: $\varepsilon\geq 0.018$.
\begin{figure*}[!t]
\centering
\includegraphics[width=160pt]{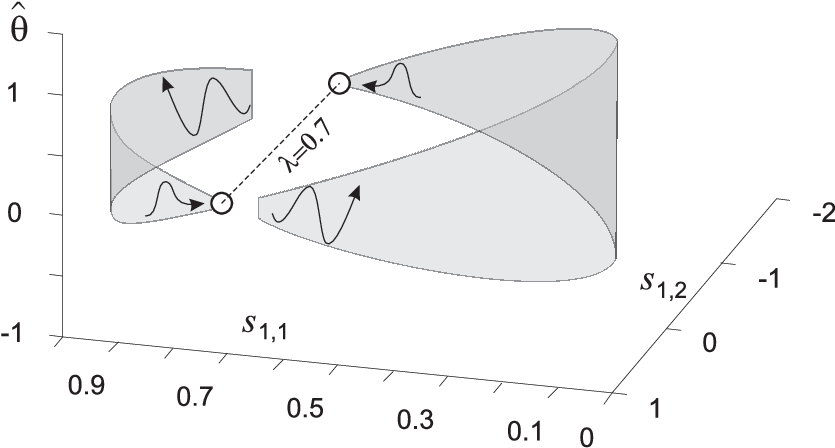}
\includegraphics[width=140pt]{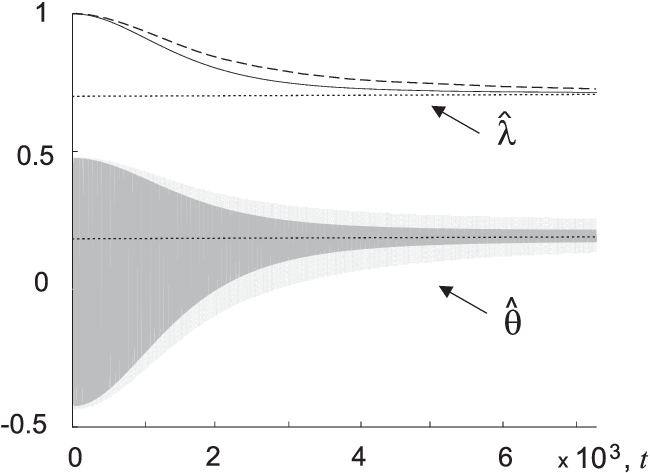}\hspace{2mm}
\includegraphics[width=160pt]{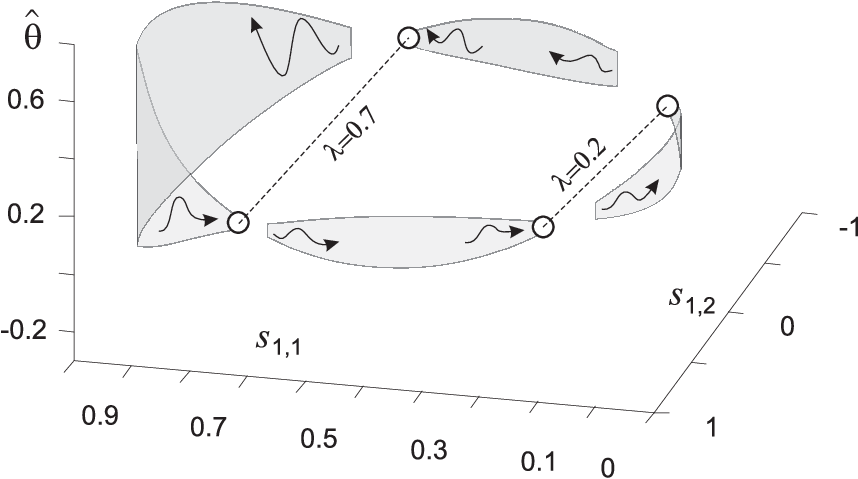}
\caption{{\it Left panel}: a qualitative picture of the system
dynamics in the coordinates $(s_{1}, s_{2},\hat{\theta})$.
Shaded regions depict envelops of $15$ trajectories  of the
system for various initial conditions. Actual trajectories of
the system are very oscillatory, and individual trajectories
are hardly distinguishable. Qualitatively, their
behavior is shown by the black arrowed lines. {\it Middle
panel}: typical simulated trajectories of $\hat{\theta}$,
$\hat{\lambda}$ as functions of $t$. Solid and dark grey curves correspond to the case when no measurement noise is added; dashed and light grey curves show trajectories of the estimates in presence of additive measurement noise. Dotted lines indicate true
values of $\theta$ and $\lambda$. {\it Right panel}: envelops
of the modified system solutions in the coordinates $(s_{1},
s_{2},\hat{\theta})$. }\label{fig:example:1}
\end{figure*}

Computer simulation results of the combined system
(\ref{eq:example_academic:1}) -- (\ref{eq:example_academic:3})
with parameters $\theta=0.2$, $\lambda=0.7$, $\gamma=0.0028$,
$\gamma_\theta=1$, $\varepsilon=0.018$ are summarized in Fig.
\ref{fig:example:1}.  As can be observed, the system has two
weakly attracting sets (marked as white circles). These sets
correspond to the true values of $\theta$ and
$\lambda$. Even though trajectories of the system are
converging to the attracting sets asymptotically, small
neighborhoods of these sets are not
forward-invariant. Hence the sets themselves are not globally asymptotically stable,
albeit they are clearly attracting. Middle panel depicts typical trajectories of $\hat{\lambda}$ and $\hat{\theta}$.
To show how the proposed observer behaves in presence of measurement noise we simulated the model-observer system in which signal $y(t)=x_1(t)$ in the observer subsystem was replaced with $y_d(t)=x_1(t)+0.05\sin(2t)$. The value of $\varepsilon$ was changed to $0.068$ to account for this perturbation. The observer retained functionality, albeit with lower precision of estimation.

In order to illustrate the
behavior of the system in the case of multiple equivalent parameterizations,
we simulated a modified version of the combined system
(\ref{eq:example_academic:1})--(\ref{eq:example_academic:3}), in
which the nonlinearly parameterized terms, i.e.
$\bfB(\sin(\lambda \cos(t))+e^{\lambda \sin(t)})$ in
(\ref{eq:example_academic:1}) and  $\bfB(\sin(\hat\lambda
\cos(t))+e^{\hat\lambda \sin(t)})$ in
(\ref{eq:example_academic:2}), are replaced with
$\bfB(\sin((\lambda-0.45)^2 \cos(t))+e^{(\lambda-0.45)^2
\sin(t)})$ and $\bfB(\sin((\hat{\lambda}-0.45)^2
\cos(t))+e^{(\hat{\lambda}-0.45)^2 \sin(t)})$ respectively.  Assumptions \ref{assume:matrices}, \ref{assume:varphi} and
A1 in Assumption \ref{assume:nonlinearities}
still hold for the modified system (with the
same values of parameters). Yet, system
(\ref{eq:example_academic:1}) with the modified
$\bfg(t,\lambda)=\bfB(\sin((\lambda-0.45)^2
\cos(t))+e^{(\lambda-0.45)^2 \sin(t)})$  is no longer  uniquely
identifiable since $\bfg(t,0.7)=\bfg(t,0.2)$ for all $t$.
Simulation results of the modified system are presented in Fig,
\ref{fig:example:1}, right panel. Instead of two attracting sets as
in the previous configuration, the modified system
has four weakly attracting sets corresponding to two equivalent
parameterizations $\theta=0.2, \lambda=0.7$ (true) and
$\theta=0.2$, $\lambda=0.2$ (spurious). Parameter estimates converge to small vicinities of
these alternative parameterizations.  Note that the estimates do not jump between neighborhoods of $\theta=0.2, \lambda=0.7$ and $\theta=0.2$, $\lambda=0.2$, which
is consistent with Remark \ref{rem:jumps}.

\begin{figure}
\centering
\includegraphics[width=0.8\columnwidth]{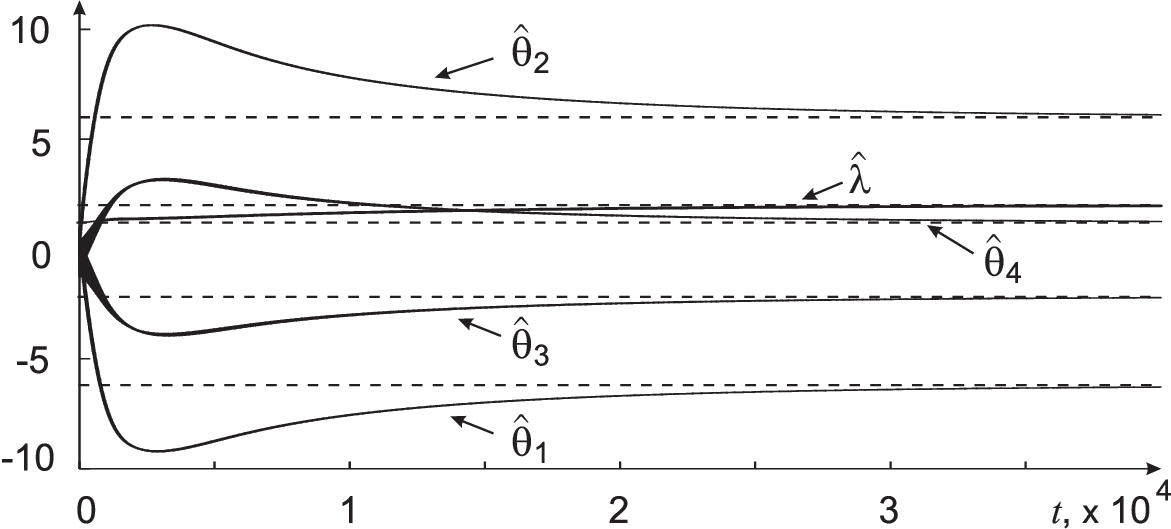}
\caption{Estimates $\hat{\thetavec}$, $\hat{\lambda}$ as functions of $t$. True values of $\thetavec$, $\lambda$
are shown as dashed lines.}\label{fig:example:RS}
\end{figure}

Finally, we illustrate
the applicability of our
approach to models (\ref{eq:neural_model:ext}). Consider the third example from Table
\ref{tab:examples} with nominal parameter
values as follows: $\tau_m=0.1666$, $\tau_s=5$, $A_f=1$, $\sigma_f=2$,
$\sigma_s=0.8$. Suppose that true values of these parameters
are unknown, but it is known that they are within $\pm 25\%$ of
their nominal values. Since the pair $\bfA,\bfC$ is observable,
there is a parameter-dependent coordinate transformation
$\bfx\mapsto T \bfx$, $T=\left(\begin{array}{cc}1 & 0 \\
\tau_s^{-1} & -\tau_m^{-1}\end{array}\right)$ rendering the
original equations into (\ref{eq:neural_model:3}) with
$\Psivec(t,\lambda,y)=\left(\begin{array}{cccc} y &
\tanh(\lambda y) & 0 & 0\\ 0 & 0 & y& \tanh(\lambda
y)\end{array}\right)$, $\bfg=0$ and
$\thetavec=\col{-\frac{1}{\tau_s}-\frac{1}{\tau_m},\frac{A_f}{\tau_m},-\frac{1+\sigma_s}{\tau_m\tau_s},\frac{A_f}{\tau_m
\tau_s}}$, $\lambda=\frac{\sigma_f}{A_f}$.  Note that $\bfA$, $\Psivec$ and $\thetavec$ differ from those in the original parametrizaton. Let
$\bfB=\col{1,1}$ and consider
$\bfM(t,[\lambda,y])=\left(m_{ij}(t,[\lambda,y])\right)$, $i=1,2$,
$j=1,\dots,4$ in (\ref{eq:linear_par_observer:modified}). It is
clear that the polynomial $s+1$ formed by the coefficients of
$\bfB$ is Hurwitz, $m_{1,j}(t,[\lambda,y])=0$,
$m_{2,j}(t,[\lambda,y])$ are defined as
$\dot{m}_{2,1}=-\dot{m}_{2,1}- y(t)$,   $\dot{m}_{2,2}=-\dot{m}_{2,2}- \tanh(\lambda y(t))$, $\dot{m}_{2,3}=-\dot{m}_{2,3}+ y(t)$, $\dot{m}_{2,4}=-\dot{m}_{2,2}+ \tanh(\lambda y(t))$,
$m_{2,j}(t_0)=0$, and that
$\varphivec(t,[\lambda,y])$ $=$ $\mathrm{col}(m_{2,1}(t,[y])+y(t), m_{2,2}(t,[\lambda,y])+\tanh(\lambda
y(t)),m_{2,3}(t,[y]), m_{2,4}(t,[\lambda,y]))$.
Given that
$A_f$, $\sigma_f$ vary within $25\%$ of their nominal values we
obtain that $\Omega_\lambda=[1.2,3.33]$. For the given system, $y$ is bounded,
 $\varphivec(t,[\lambda,y])$, as a function of $t,\lambdavec$ on $\Real_{\geq t_0}\times\Omega_\lambda$, is $\lambda$-UPE with
$T=100$, $\mu=0.08$. Moreover the restriction of $\alphavec_5$,
defined in (\ref{eq:nonlinear_pe_new}), on $\Real_{\geq t_0}\times\Omega_\lambda$ is $\lambda'$-UPE with
$T=100$, $\mu=0.0054$. Hence assumptions of Theorem
\ref{cor:general} are met. We simulated the system and observer
(\ref{eq:linear_par_observer:modified}),
(\ref{eq:nonlinear_par_observer}) with $\gamma_{\theta}=4$,
$\gamma=0.004$, and $\varepsilon=0$  for various initial
conditions and values of $\thetavec$, $\lambda$;
$\hat{\thetavec}$, $\hat{\lambda}$ approached true values of
$\thetavec$, $\lambda$ asymptotically as prescribed. An example
of typical behavior of the estimates is shown in Fig.
\ref{fig:example:RS}. Original parameters of the model can be
recovered as: $\hat{\tau}_s=\hat{\theta}_2/\hat{\theta}_4$,
$\hat{\tau}_m=-1/(\hat{\theta}_1+1/\hat{\tau}_s)$,
$\hat{A}_f=\hat{\tau}_m \hat{\theta}_2$,
$\hat{\sigma}_s=-\hat{\theta}_3\hat{\tau}_s\hat{\tau}_m-1$,
$\hat{\sigma}_f=\hat{A}_f\hat{\lambda}$.
Further examples may be found in the supplementary material \cite{Tyukin:Supplement:2013}.


\section{Conclusion}\label{sec:conclusion}

We derived an observer that can reconstruct asymptotically the
unknown state and parameter values of  a class of systems with
general nonlinear parametrization. This class can be viewed as
an extension of the adaptive observer canonical forms
\cite{Bastin88}, \cite{Marino90}. In contrast to earlier
approaches addressing the problem of nonlinear parametrization
in the problem of adaptive observer design
\cite{Automatica:Farza:2009},\cite{Automatica:Grip:2010},
\cite{IEEE:2009:Grip}, \cite{SysConLett:Grip:2011},
\cite{Ortega2011}, \cite{IEEE_TAC:2007:tyukin}, the class of
parameterizations for which the reconstruction is
guaranteed is not limited to convex/concave or one-to-one
functions.  We showed that reconstruction can be
achieved, subject to additional conditions of linear/nonlinear
persistency of excitation, if nonlinearly parameterized
functions in the model are bounded, differentiable and
Lipschitz.

The set to which the estimates
converge is not guaranteed to be asymptotically stable. Yet the
set is attracting in a weaker, Milnor sense, cf.
\cite{Milnor_1985}. Numerical simulations revealed that the convergence time in our
approach depends heavily on the dimension of
$\lambdavec$; it does not, however, depend crucially on the
dimension of $\thetavec$. 
This renders the method more efficient than exhaustive search;
the smaller the dimension of $\lambdavec$
the more advantageous our method becomes.  In this
respect a related question arises: is there a ``best''
parametrization  for a given physical model in the class of
systems (\ref{eq:neural_model}) or (\ref{eq:neural_model:ext})?
The answer is likely to require quantitative assessment of the
performance of various observers for all admissible
parametrizations. We do not answer this question here, but hope
to be able to address it in future.


\bibliographystyle{plain}
\bibliography{adaptive_observer_automatica}


\appendix
\section{Appendix}\label{app:proofs}

\begin{lem}\label{lem:system_io_inferrence_1order} Consider $\dot{y}=k y + u(t)+d(t)$, $k\in\Real$, $u,d:\Real_{\geq t_0}\rightarrow\Real$, $u\in\mathcal{C}^1$, $ d\in\mathcal{C}^{0}$, and let $\max\{|u(t)|,|\dot{u}(t)|\}\leq B$, $|d(t)|\leq \Delta_\xi$. Then
$\|y\|_{\infty,[t_0,\infty)}\leq \varepsilon \Rightarrow  \
\exists  \ t_1(\varepsilon)\geq t_0: \
\|u\|_{\infty,[t_1(\varepsilon),\infty)} \leq
\sqrt{\varepsilon}(1+e^{|k|\sqrt{\varepsilon}}+B)+\Delta_\xi$.
\end{lem}
\begin{proof}[Proof of Lemma \ref{lem:system_io_inferrence_1order}.]
Noticing that $y(t)$ for $t\geq t_0 + T$, $T>0$, can be expressed as:
$y(t)=y(t-T)e^{ k T}+\int_{t-T}^t
e^{k(t-\tau)}(u(\tau)+d(\tau))d\tau$ and using the Mean-value
theorem we obtain: $y(t)-y(t-T)e^{ k T}= T
e^{k(t-\tau')}(u(\tau')+d(\tau')), \ \tau'\in[t-T,t]$. Hence  $ \varepsilon(1+e^{kT})\geq T e^{k(t-\tau')}
(|u(t)|-T B-\Delta_{\xi})$, and
$\Delta_\xi+TB+\frac{\varepsilon(1+e^{kT})}{T \min\{1,e^{kT}\}} \geq
\Delta_\xi+TB+\frac{\varepsilon(1+e^{kT})}{T \min\{1,e^{k(t-\tau')}\}}
\geq |u(t)|$   for all  $t\geq t_0+T$.
Given that $T$ can be chosen arbitrarily we let
$T=\sqrt{\varepsilon}$, and thus $|u(t)|\leq
\sqrt{\varepsilon}(1+e^{k\sqrt{\varepsilon}})\max\{1,e^{-k\sqrt{\varepsilon}}\}+B\sqrt{\varepsilon}
+\Delta_\xi \leq
\sqrt{\varepsilon}(1+e^{|k|\sqrt{\varepsilon}}+B)+\Delta_{\xi}
\ \forall \ t\geq t_0 + \sqrt{\varepsilon}$. $\square$
\end{proof}

\begin{proof}[Proof of Lemma \ref{lem:observer_inferrence}.]
Let us rewrite (\ref{eq:system_io}) as
\[
\begin{split}
\dot{y}&= a_{1} y + \tilde{\bfC}\tilde\bfx + u_1(t)+d_1(t)\\
\dot{\tilde \bfx}&= \tilde\bfA \tilde \bfx + \tilde \bfa y +\bfb u_1(t) + \bfG \bfu(t)+\tilde{\bfd}(t), \end{split}
\]
where $\tilde{\bfa}=\mathrm{col}(a_2,\dots,a_n)$,
$\tilde\bfC=\mathrm{col}(1,0,\dots,0)$,
$\tilde\bfd(t)=\mathrm{col}(d_2(t),\dots,d_n(t))$, and
$\bfG=\left(\begin{array}{cc} - {\bfb} & I_{n-1}
\end{array}\right)$,  $\tilde  \bfA=\left(\begin{array}{cc}0 & I_{n-2}\\
                                0 & 0
                                \end{array}\right)$.
Let $\|y(t)\|_{\infty,[t_0,\infty)}\leq \varepsilon$ and denote
$e(t)=\tilde\bfC^{T}\tilde \bfx(t)+u_1(t)$. According to Lemma
\ref{lem:system_io_inferrence_1order},  there is a
$t_1(\varepsilon) > t_0$  and
$\upsilon_1,\upsilon_2\in\mathcal{K}$ such that
$\|e(t)\|=\|\tilde\bfC^{T}\tilde\bfx(t)+u_1(t)\|\leq
\upsilon_1(\varepsilon)+\upsilon_2(\Delta_\xi) \ \forall  \
t\geq t_1(\varepsilon)$.

 Using the notation above we obtain:
$\dot{\tilde\bfx}=(\tilde \bfA -  \bfb
\tilde\bfC^{T})\tilde\bfx + \tilde \bfa y(t) +
\tilde{\bfG}\bfu(t) +  \bfb e(t)+\tilde{\bfd}(t)$.
Matrix $\tilde \bfA -  \bfb \tilde\bfC^{T}=\Lambda$ is Hurwitz,
and hence there are $D,k\in\Real_{>0}$ such that
$\|e^{\Lambda(t-t_0)}\|\leq D e^{-k (t-t_0)}$. {
Therefore $ \|\tilde \bfC^{T} \tilde\bfx(t)-
\tilde\bfC^{T}\int_{t_0}^{t}e^{\Lambda(t-\tau)}\bfG
\bfu(\tau)d\tau\| \leq De^{-k(t-t_0)}\|\tilde{\bfx}(t_0)\|
+\frac{D}{k}(\|\bfa\|\varepsilon
+\|\bfb\|(\upsilon_1(\varepsilon)+\upsilon_2(\Delta_\xi))+\Delta_\xi)
$.

Noticing that
$z_1=\tilde\bfC^{T}\int_{t_0}^{t}e^{\Lambda(t-\tau)}\bfG
\bfu(\tau)d\tau$ and denoting $\kappa_1(\varepsilon)=2
\frac{D}{k}(\|\bfa\|\varepsilon
+\|\bfb\|\upsilon_1(\varepsilon))+\upsilon_1(\varepsilon)$,
$\kappa_2(\Delta_\xi)=2
\frac{D}{k}(\Delta_\xi+\|\bfb\|\upsilon_2(\Delta_\xi))+\upsilon_2(\Delta_\xi)$
} we can conclude that there is a $t'(\varepsilon,\bfx_0)\geq
t_1(\varepsilon)$ such that
\[
\|z_1(\tau)+ u_1(\tau)\|_{\infty,[t,\infty)}\leq \kappa_1(\varepsilon)+\kappa_2(\Delta_\xi) \ \forall \ t\geq t'(\varepsilon).
\]
Noticing that $y(t),\bfd(t)\equiv 0 \Rightarrow e(t)\equiv0$ ensures that (\ref{eq:identifiability}) holds too. \hfill $\square$
\end{proof}

\begin{proof}[Proof of Lemma \ref{lem:UPE}.]
Consider
$\bfJ(\lambdavec,t)$ $=$ $\bfz^{T}
\big(\int_t^{t+T}$ $\alphavec(\tau,\lambdavec)\alphavec^T(\tau,\lambdavec){\rm
d}\tau \big)\bfz=\int_{t}^{t+T}
\|\bfz^{T}\alphavec(\tau,\lambdavec)\|^2$,  where $\bfz\in\Real^{n+m}$, $\bfz\neq 0$, for $t\in\Real_{\geq t_0}$. According to C2 we have: $\bfJ(\lambdavec,t)\geq \mu
\|\bfz\|^2 \ \forall \ \lambdavec\in\Omega_{\lambda}$. Let  $\hat \lambdavec:\Real_{\geq t_0}\rightarrow\Omega_\lambda$ be a
differentiable function, and consider
$\bfJ(\hat\lambdavec(t),t) -
\int_t^{t+T}\|\bfz^{T}\alphavec(\tau,\hat\lambdavec(\tau))\|^2{\rm
d}\tau$
 $=$ $\int_t^{t+T}\|\bfz^{T}\alphavec(\tau,\hat\lambdavec(t))\|^2$
 $-$ $\|\bfz^{T}\alphavec(\tau,\hat\lambdavec(\tau))\|^2 {\rm d}\tau$
 $=$ $\int_t^{t+T}\|\bfz^{T}\alphavec(\tau,\hat\lambdavec(t))\|^2 - \bfz^{T} \alphavec(\tau,\hat\lambdavec(t)) \alphavec^{T}(\tau,\hat\lambdavec(\tau))\bfz$
  $+$
  $\bfz^{T}\alphavec(\tau,\hat\lambdavec(t)) \alphavec^{T}(\tau,\hat\lambdavec(\tau))\bfz- \|\bfz^{T}\alphavec(\tau,\hat\lambdavec(\tau))\|^2 {\rm d}\tau$
  $=$ $\int_t^{t+T}\bfz^{T}\alphavec(\tau,\hat\lambdavec(t))[\alphavec^{T}(\tau,\hat\lambdavec(t))-\alphavec^{T}(\tau,\hat\lambdavec(\tau))]\bfz$
  $+$ $\int_t^{t+T}\bfz^{T}[\alphavec(\tau,\hat\lambdavec(t))$ $-$ $\alphavec(\tau,\hat\lambdavec(\tau))]$ $\alphavec^{T}(\tau,\hat\lambdavec(\tau))\bfz
 {\rm d}\tau$.
Applying  the Cauchy-Schwarz inequality to the last equality,
and invoking C4 and C3 we obtain: $\bfJ(\hat\lambdavec(t),t) -
\int_t^{t+T}\|\bfz^{T}\alphavec(\tau,\hat\lambdavec(\tau))\|^2{\rm
d}\tau\leq$  $
\big(\int_t^{t+T}\|\bfz^{T}[\alphavec(\tau,\hat\lambdavec(t))-\alphavec(\tau,\hat\lambdavec(\tau))]\|^2d\tau\big)^{\frac{1}{2}}$
$2 M T \|\bfz\|$ $\leq$ $2 D M T^2$ $\|\bfz\|^2$ $\hspace{5mm}$
$\max_{\tau\in[t,t+T]}\|\dot{\hat \lambdavec}(\tau)\|$. Thus
(\ref{eq:lambda_condition}), (\ref{eq:gamma_boundary:1}) ensure
that
$\int_t^{t+T}\|\bfz^{T}\alphavec(\tau,\hat\lambdavec(\tau))\|^2{\rm
d}\tau$ $\geq$ $\bfJ(\hat\lambdavec(t),t)-r \mu\|\bfz\|^2$. This, in accordance with C2, guarantees that
$\int_t^{t+T}\|\bfz^{T}\alphavec(\tau,\hat\lambdavec(\tau))\|^2{\rm
d}\tau\geq (1-r)\mu \|\bfz\|^2$. Hence $\alphavec(t,\hat\lambdavec(t))$ is persistently
exciting in the sense of Definition \ref{defn:pe}. The value of $(1-r)\mu$ does not depend on the choice of $\hat
\lambdavec$ as long as (\ref{eq:gamma_boundary:1}),
(\ref{eq:lambda_condition}) hold. Finally, notice that C4
and (\ref{eq:lambda_condition}) guarantee boundedness
of $\alphavec(\cdot,\hat\lambdavec(\cdot))$ and its derivative:
$\max\{\|\alphavec(t,\hat\lambdavec(t))\|,\|\dot{\alphavec}(t,\hat
\lambdavec(t))\|\} \leq M + M M_\lambda = M + \frac{\mu r}{2 D
T^2}$. Taking  C1  and Theorem
\ref{thm:preliminary:exponential_skew_symmetric} into account we
conclude that the lemma  follows. \hfill $\square$
\end{proof}

\begin{proof}[Proof of Lemma \ref{lem:bound}.]
According to conditions of the lemma,  (\ref{eq:gamma_choice}), we conclude that $h(t_0)\geq 0$. Introduce a strictly
decreasing sequence: $ \{\sigma_i\},  \ i=0,1,\dots,  \
\sigma_i=(1/\kappa)^i, \ \kappa\in(1,\infty)$. Further, let $
\{t_i\}, \ i=1,\dots, \ t_1<t_2<\cdots<t_n<\cdots$ be an
ordered infinite sequence:
\begin{equation}\label{eq:h_i}
 h(t_i)=\sigma_i h(t_0).
\end{equation}
If $\{t_i\}$ satisfying (\ref{eq:h_i}) does not exist then it is clear that $h(t_0)\geq h(t)\geq 0$ for all $t\geq
t_0$. Hence, in
accordance with (\ref{eq:initial_conditions}),
$\bfx(\cdot)$ is bounded for all $t\geq t_0$, and nothing remains
to be proven. Let us now show that if (\ref{eq:epsilon_choice}),
(\ref{eq:gamma_choice}), and (\ref{eq:h_i}) hold then
\begin{equation}\label{eq:limit_h}
h(t)\rightarrow 0 \Rightarrow t\rightarrow\infty.
\end{equation}
Consider $T_i=t_i-t_{i-1}$. It is clear from (\ref{eq:gamma_0}) that
\begin{equation}\label{eq:T_i_defn}
T_i D_\gamma \max_{\tau\in[t_{i-1},t_i]}\|\bfx(\tau)+\bfd(\tau)\|_{\varepsilon}\geq h(t_0)(\sigma_{i-1}-\sigma_{i}).
\end{equation}
In addition,
$\max_{\tau\in[t_{i-1},t_i]}$ $\|\bfx(\tau) $+ $\bfd(\tau)\|_{\varepsilon}$ $=$ $\|\bfx(\tau)$ $+$ $\bfd(\tau)\|_{\infty,[t_{i-1},t_i]}$ $-\varepsilon$
if
$\|\bfx(\tau)+\bfd(\tau)\|_{\infty,[t_{i-1},t_i]}>\varepsilon$,
and
$\max_{\tau\in[t_{i-1},t_i]}\|\bfx(\tau)+\bfd(\tau)\|_{\varepsilon}=0$
overwise, we can see from (\ref{eq:T_i_defn}) that
\begin{equation}\label{eq:T_i_defn:1}
T_i\geq \left\{ \begin{array}{l} \frac{h(t_0)(\sigma_{i-1}-\sigma_i)}{D_\gamma}\frac{1}{\|\bfx(\tau)+\bfd(\tau)\|_{\infty,[t_{i-1},t_i]}-\varepsilon},\\ \ \  \ \ \ \  \|\bfx(\tau)+\bfd(\tau)\|_{\infty,[t_{i-1},t_i]}>\varepsilon;\\
\infty,  \   \|\bfx(\tau)+\bfd(\tau)\|_{\infty,[t_{i-1},t_i]}\leq\varepsilon.
\end{array}\right.
\end{equation}
Consider the case when
$\|\bfx(\tau)+\bfd(\tau)\|_{\infty,[t_{i-1},t_i]}-\varepsilon>0$
for all $i=1,2,\dots$. Let us pick
\begin{equation}\label{eq:tau_choice}
 \tau^\ast={\varrho}^{-1}\left({d}/{\kappa}\right), \ d\in(0,1),
\end{equation}
and select the value of $D_\gamma$ such that
(\ref{eq:gamma_choice}) holds. Given that
$\|\bfx(\tau)+\bfd(\tau)\|_{\infty,[t_{0},t_1]}-\varepsilon\leq
{\varrho}(0)\|\bfx(t_0)\|+c h(t_0)+\Delta+\Delta_d-\varepsilon$,
conditions (\ref{eq:gamma_choice}), (\ref{eq:epsilon_choice}),
and (\ref{eq:tau_choice}) imply $D_\gamma\leq
\frac{\kappa-1}{\kappa}\frac{h(t_0)}{{\varrho}(0)\|\bfx(t_0)\|+c
|h(t_0)|}\frac{1}{\tau^\ast} \leq
\frac{h(t_0)(\sigma_0-\sigma_1)}{(\|\bfx(\tau)+\bfd(\tau)\|_{\infty,[t_{0},t_1]}-\varepsilon)}\frac{1}{\tau^\ast}$.
This, as follows from (\ref{eq:T_i_defn:1}), guarantees that
$T_1\geq \tau^\ast$.

Without loss of generality suppose that there is  an $i\geq 2$: $T_j\geq \tau^\ast$ for all $1\leq j\leq i-1$. Let us show that
$T_{i-1}\geq \tau^\ast\Rightarrow T_{i}\geq \tau^\ast$. This
will ensure that (\ref{eq:limit_h}) is satisfied and that the lemma hold. Consider
$\|\bfx(\tau)\|_{\infty,[t_{i-1},t_i]}$;
(\ref{eq:interconnection}) and (\ref{eq:h_i}) imply that:
$\|\bfx(\tau)\|_{\infty,[t_{i-1},t_i]}\leq
{\varrho}(0)\|\bfx(t_{i-1})\|+c\sigma_{i-1}h(t_0)+\Delta$. Hence
\[
\begin{split}
&\|\bfx(\tau)\|_{\infty,[t_{i-1},t_i]}\leq{\varrho}(0)[{\varrho}(T_{i-1})\|\bfx(t_{i-2})\|+c\sigma_{i-2}h(t_0)]\\
&+{\varrho}(0)\Delta +ch(t_0)\sigma_{i-1}+\Delta\leq { \varrho}(0){ \varrho}(\tau^\ast) \|\bfx(t_{i-2})\| + P_1,
\end{split}
\]
where $P_1={ \varrho}(0)c
\sigma_{i-2}h(t_0)+c\sigma_{i-1}h(t_0)+{ \varrho}(0)\Delta+\Delta$.
Invoking (\ref{eq:interconnection}) again results in
\[
\begin{split}
&\|\bfx(\tau)\|_{\infty,[t_{i-1},t_i]}\leq{ \varrho}(0){ \varrho}^2(\tau^\ast)\|\bfx(t_{i-3})\|+P_2,
\end{split}
\]
where $P_2=c
h(t_0){ \varrho}(0)[{ \varrho}(\tau^\ast)\sigma_{i-3}+\sigma_{i-2}]+c\sigma_{i-1}h(t_0)+
 { \varrho}(0)[{ \varrho}(\tau^\ast)\Delta+\Delta]+\Delta$, and
\[
\begin{split}
&\|\bfx(\tau)\|_{\infty,[t_{i-1},t_i]}\leq{ \varrho}(0){ \varrho}^3(\tau^\ast)\|\bfx(t_{i-4})\|+P_3,
\end{split}
\]
where $P_3=c h(t_0){ \varrho}(0)[{ \varrho}^2(\tau^\ast)\sigma_{i-4}+{ \varrho}(\tau^\ast)\sigma_{i-3}+\sigma_{i-2}]+c\sigma_{i-1}h(t_0)+ \Delta{ \varrho}(0)[{ \varrho}(\tau^\ast)^2+{ \varrho}(\tau^\ast)+1]+\Delta=c
h(t_0){ \varrho}(0)[\sum_{j=0}^2{ \varrho}^j(\tau^\ast)\sigma_{i-j-2}]+ch(t_0)\sigma_{i-1}+\Delta
{ \varrho}(0)[\sum_{j=0}^2{ \varrho}^j(\tau^\ast)]+\Delta$.
After $i-1$ steps we obtain
\begin{equation}\label{eq:x_i_est}
\begin{split}
\|\bfx(\tau)\|_{\infty,[t_{i-1},t_i]}&\leq { \varrho}(0){ \varrho}^{i-1}(\tau^\ast)\|\bfx(t_0)\|+P_{i-1},
\end{split}
\end{equation}
with
$P_{i-1}=ch(t_0){ \varrho}(0)[\sum_{j=0}^{i-2}{ \varrho}^j(\tau^\ast)\sigma_{i-j-2}]+ch(t_0)\sigma_{i-1}$ $+$ $\Delta
{ \varrho}(0)[\sum_{j=0}^{i-2}{ \varrho}^j(\tau^\ast)]+\Delta$.
The values of $T_i$, as follows from (\ref{eq:T_i_defn:1}), are
bounded from below by:
\begin{equation}\label{eq:T_i_bound}
\begin{array}{l}
T_i\geq \frac{\sigma_{i-1}-\sigma_{i}}{\sigma_{i-1} D_{\gamma}\sigma_{i-1}^{-1}}\frac{h(t_0)}{\left(\|\bfx(\tau)+\bfd(\tau)\|_{\infty,[t_{i-1},t_i]}-\varepsilon\right)}.
\end{array}
\end{equation}
Consider
$\sigma_{i-1}^{-1}(\|\bfx(\tau)+\bfd(\tau)\|_{\infty,[t_{i-1},t_i]}-\varepsilon)$.
Taking (\ref{eq:x_i_est}) into account we derive that:
\[
\begin{array}{l}
\sigma_{i-1}^{-1}(\|\bfx(\tau)+\bfd(\tau)\|_{\infty,[t_{i-1},t_i]}-\varepsilon)\leq { \varrho}(0) { \varrho}^{i-1}(\tau^\ast)\times\\
\kappa^{i-1}\|\bfx(t_0)\|+k^{i-1}P_{i-1}+k^{i-1}\Delta_d-k^{i-1}\varepsilon= \\
{ \varrho}(0) { \varrho}^{i-1}(\tau^\ast)\kappa^{i-1}\|\bfx(t_0)\|+ch(t_0){ \varrho}(0)\kappa[\sum_{j=0}^{i-2}{ \varrho}^j(\tau^\ast)\kappa^j]\\
+ch(t_0)+\kappa^{i-1}[\Delta({ \varrho}(0)\sum_{j=0}^{i-2}{ \varrho}^j(\tau^\ast)+1)+\Delta_d-\varepsilon].
\end{array}
\]
Noticing that $\tau^\ast$ is chosen in accordance with
(\ref{eq:tau_choice}) one can therefore obtain:
\[
\begin{array}{l}
\sigma_{i-1}^{-1}(\|\bfx(\tau)+\bfd(\tau)\|_{\infty,[t_{i-1},t_i]}-\varepsilon)\leq { \varrho}(0)\|\bfx(t_0)\|+\\
ch(t_0)+ch(t_0){ \varrho}(0)\kappa\sum_{j=0}^{i-2}d^j+\\
\kappa^{i-1}[\Delta({ \varrho}(0)\sum_{j=0}^{i-2}\frac{d}{\kappa}^j+1)+\Delta_d-\varepsilon]\\
\leq { \varrho}(0)\|\bfx(t_0)\|+ch(t_0)(1+\frac{{ \varrho}(0)\kappa}{1-d})+\\
k^{i-1}[\Delta(\frac{{ \varrho}(0)}{1-d/k}+1)+\Delta_d-\varepsilon].
\end{array}
\]
Condition (\ref{eq:epsilon_choice}) implies that
$\Delta(\frac{{ \varrho}(0)}{1-d/k}+1)+\Delta_d-\varepsilon\leq
0$. Hence
$\sigma_{i-1}^{-1}(\|\bfx(\tau)+\bfd(\tau)\|_{\infty,[t_{i-1},t_i]}-\varepsilon)\leq
{ \varrho}(0)\|\bfx(t_0)\|+ch(t_0)(1+\frac{{ \varrho}(0)\kappa}{1-d})$. Substituting the latter estimate into
(\ref{eq:T_i_bound}) and using (\ref{eq:gamma_choice}) yields
$T_{i}\geq  \tau^\ast$.
Thus $h(\cdot)$ is bounded for $t\geq t_0$, and
hence so is $\bfx(\cdot)$. \hfill $\square$
\end{proof}

\begin{proof}[Proof of Theorem \ref{cor:general}.]
Let $\Lambda_0=\bfA-\bfB\bfC^{T}\bfA$,
$\bfG_0=\bfI_n-\bfB\bfC^{T}$. Consider
$\tilde{\varphivec}(t,\hat{\lambdavec}(t),T_1)$ $=$
$\bfC^{T}\bfA\int_{t-T_1}^{t}$ $e^{\Lambda_0(t-\tau)}\bfG_0$ $\Psivec(\tau,\hat{\lambdavec}(t),y(\tau))d\tau$ $+$ $\bfC^{T}\Psivec(t,\hat{\lambdavec}(t),y(t))$.
It is clear that for any $\varepsilon_1>0$ there are $T_1$,
$t_1$ sufficiently large and $\gamma_1$ sufficiently small:
\begin{equation}\label{eq:proofs:filter_PE}
\|\tilde{\varphivec}(t,\hat{\lambdavec}(t),T_1)-\varphivec(t,\hat{\lambdavec}(t),y(t),[\hat\lambdavec,y])\|<\varepsilon_1
\end{equation}
for all $\gamma\in(0,\gamma_1)$ and $t\geq t_1\geq t_0$. Indeed, consider
$\delta_1(T_1,t)=\bfC^{T}\bfA e^{\Lambda_0 T_1}\int_{t_0}^{t-T_1}$ $e^{\Lambda_0(t-T_1-\tau)}\bfG_0$ $\Psivec(\tau,\hat{\lambdavec}(\tau),y(\tau))d\tau$,
$\delta_2(T_1,t)=\bfC^{T}\bfA\int_{t-T_1}^{t} e^{\Lambda_0(t-\tau)}$ $\bfG_0
(\Psivec(\tau,\hat{\lambdavec}(\tau),y(\tau))$ $-$ $\Psivec(\tau,\hat{\lambdavec}(t),y(\tau)))d\tau$, pick $\varepsilon_1>0$, and let  $T_1$ be so large that $|\delta_1(T_1,t)|< \varepsilon_1 / 2$ for $t\geq t_0+T_1$. Let $\gamma_1\in\Real_{>0}$ be so small that $|\delta_2(T_1,t)|<\varepsilon_1/2$ for all $t\geq t_0+T_1$ (such a choice is always possible  due to that $\Psivec$ is Lipschitz in $\hat\lambdavec$). Noticing that $\varphivec(t,\hat{\lambdavec}(t),y(t),[\hat\lambdavec,y])=\delta_1(T_1,t)+\delta_2(T_1,t)+\tilde{\varphivec}(t,\hat{\lambdavec}(t),T_1)$ we can conclude that (\ref{eq:proofs:filter_PE}) holds.

Given that the restriction of
$\varphivec(t,\lambdavec,y(t),[\lambdavec,y])$ (as a function of $t$, $\lambdavec$) on $\Real_{\geq t_0}\times\Omega_\lambda$ is $\lambda$-UPE with $T,\mu$,
there is a $\varepsilon_1$ in (\ref{eq:proofs:filter_PE}) such that
$\tilde{\varphivec}(t,{\lambdavec},T_1)\in\lambda\mathrm{UPE}(T,\mu-\epsilon)$,
$\epsilon\in(0,\mu/3)$. On the other hand (see the first part
of the proof of Lemma \ref{lem:UPE}), there is a $\gamma_2$
such that $\tilde{\varphivec}(t,\hat{\lambdavec}(t),T_1)$ is
persistently exciting with parameters $T,\mu-2\epsilon$ for all
$\gamma\in(0,\gamma_2)$. Choosing
$\gamma\in(0,\min\{\gamma_1,\gamma_2\})$ and taking
(\ref{eq:proofs:filter_PE}) into account we conclude that the restriction of
$\varphivec(\cdot,\hat{\lambdavec}(\cdot),y(\cdot),[\hat\lambdavec,y])$ on $\Real_{\geq t_2}$ is
persistently exciting ($t_2>t_1>t_0$) provided that
$\varepsilon_1$ is small enough and $t_2$ is sufficiently
large.  Thus, invoking the argument presented in Part 2 of the
proof of Theorem \ref{theorem:neural_identification} we can
conclude that (\ref{eq:y_limit}) and
(\ref{eq:convervence_lambda_h}) hold for the combined system.
Convergence of state and parameter estimates can now be shown similarly to the 3d part of the proof
of Theorem \ref{theorem:neural_identification}. \hfill $\square$
\end{proof}

\end{document}